\documentclass[11pt]{article}
\usepackage{amsmath}
\usepackage{amsfonts}
\usepackage{graphicx}
\usepackage{amssymb}
\usepackage{epsfig}
\usepackage{color,xcolor}

\def\sqr#1#2{{\vcenter{\vbox{\hrule height.#2pt
              \hbox{\vrule width.#2pt height#1pt \kern#1pt \vrule width.#2pt}
          \hrule height.#2pt}}}}
\def\signed #1{{\unskip\nobreak\hfil\penalty50
          \hskip2em\hbox{}\nobreak\hfil#1
          \parfillskip=0pt \finalhyphendemerits=0 \par}}
\def\endpf{\signed {$\sqr69$}}
\def\sqr#1#2{{\vcenter{\vbox{\hrule height.#2pt
              \hbox{\vrule width.#2pt height#1pt \kern#1pt \vrule width.#2pt}
              \hrule height.#2pt}}}}
\def\signed #1{{\unskip\nobreak\hfil\penalty50
              \hskip2em\hbox{}\nobreak\hfil#1
              \parfillskip=0pt \finalhyphendemerits=0 \par}}
\def\endpf{\signed {$\sqr69$}}
\def\3n{\negthinspace \negthinspace \negthinspace }
\def\2n{\negthinspace \negthinspace }
\def\1n{\negthinspace }

\def\={\buildrel \triangle \over =}

\def\O{\Omega}

\def\bs{\bigskip}
\def\q{\quad}
\def\qq{\qquad}

\def\max{\mathop{\rm max}}
\def\min{\mathop{\rm min}}

\def\sup{\mathop{\rm sup}}

\def\inf{\hbox{\rm inf$\,$}}

\def\|{\Big |}
\def\({\Big (}
\def\){\Big )}
\def\[{\Big[}
\def\]{\Big]}
\def\be{\begin{equation}}
\def\bel{\begin{equation}\label}
\def\ee{\end{equation}}
\def\bt{\begin{theorem}}
\def\bcd{\begin{condition}}
\def\ecd{\end{condition}}
\def\et{\end{theorem}}
\def\bc{\begin{corollary}}
\def\ec{\end{corollary}}
\def\bde{\begin{definition}}
\def\ede{\end{definition}}
\def\bl{\begin{lemma}}
\def\el{\end{lemma}}
\def\bp{\begin{proposition}}
\def\ep{\end{proposition}}
\def\br{\begin{remark}}
\def\er{\end{remark}}
\def\ba{\begin{array}}
\def\ea{\end{array}}
\def\ed{\end{document}}

\def\square#1{\vbox{\hrule\hbox{\vrule height#1%
     \kern#1\vrule}\hrule}}
\def\rectangle#1#2{\vbox{\hrule\hbox{\vrule height#1%
     \kern#2\vrule}\hrule}}

\font\tenbb=msbm10 \font\sevenbb=msbm7 \font\fivebb=msbm5

\newfam\bbfam
\scriptscriptfont\bbfam=\fivebb \textfont\bbfam=\tenbb
\scriptfont\bbfam=\sevenbb

\newtheorem{lemma}{Lemma}[section]
\newtheorem{remark}{Remark}[section]
\newtheorem{example}{Example}[section]
\newtheorem{theorem}{Theorem}[section]
\newtheorem{corollary}{Corollary}[section]

\newtheorem{definition}{Definition}[section]
\newtheorem{proposition}{Proposition}[section]
\newtheorem{condition}{Condition}[section]

\makeatletter
   
   \@addtoreset{equation}{section}
\makeatother

\begin{document}
\title{Regularity properties for general HJB equations. A BSDE method}

\author{Rainer Buckdahn{\footnote{This work has been done in the frame of the Marie Curie ITN Project ``Deterministic and Stochastic Controlled Systems and Applications'', call: F97-PEOPLE-2007-1-1-ITN, no: 213841-2. \qq \qq \qq \qq \qq  \mbox{ }  \mbox{ } \mbox{ }  \mbox{ } \mbox{ }  \mbox{ }\mbox{ }  \mbox{ }\mbox{ }  \mbox{ } \mbox{   } 
$^{\ast\ast}$Jianhui Huang acknowledges the financial support from the RGC Earmarked Grants Poly U. 500909 and 501010.\qq \qq \qq  \mbox{ }
$^{\ast\ast\ast}$Juan Li is the corresponding author, and has been supported by the NSF of P.R.China (No. 11071144), Shandong
Province (No. BS2011SF010), Independent Innovation Foundation of Shandong University, and National Basic Research Program of China (973 Program) (No. 2007CB814904), 111
Project (No. B12023).}}\\
{\small D\'{e}partement de Math\'{e}matiques, Universit\'{e} de
Bretagne Occidentale,}\\
 {\small 6, avenue Victor-le-Gorgeu, CS 93837, 29238 Brest
Cedex 3, France.}\\
{\small{\it E-mail: Rainer.Buckdahn@univ-brest.fr.}}\\
Jianhui Huang$^{\ast\ast}$\\
{\small Department of Applied Mathematics,}\\
{\small The Hong Kong Polytechnic University, P. R. China.}\\
{\small{\it E-mail: majhuang@inet.polyu.edu.hk.}}\\
 Juan Li$^{\ast\ast\ast}$\\
{\small School of Mathematics and Statistics,}\\
{\small Shandong University at Weihai, Weihai 264209, P. R. China.}\\
{\small {\it E-mail: juanli@sdu.edu.cn.}} \date{April 23, 2011}}
\maketitle

\abstract{In this work we investigate regularity properties of a large class of Hamilton-Jacobi-Bellman (HJB) equations with or without obstacles, which can be stochastically interpreted in form of a stochastic control system which nonlinear cost functional is defined with the help of a backward stochastic differential equation (BSDE) or a reflected BSDE (RBSDE). More precisely, we prove that, firstly, the unique viscosity solution $V(t,x)$ of such a HJB equation  over the time interval $[0,T],$ with or without an obstacle, and with terminal condition at time $T$, is jointly Lipschitz in $(t,x)$, for $t$ running any compact subinterval of $[0,T)$. Secondly, for the case that $V$ solves a HJB equation without an obstacle or with an upper obstacle it is shown under appropriate assumptions that $V(t,x)$ is jointly semiconcave in $(t,x)$. These results extend earlier ones by Buckdahn, Cannarsa and Quincampoix [1]. Our approach embeds their idea of time change into a BSDE analysis. We also provide an elementary counter-example which shows that, in general, for the case that $V$ solves a HJB equation with a lower obstacle the semi-concavity doesn't hold true.}

 \medskip
 \noindent{{\bf AMS subject classification.} 93E20, 35D40, 60H10, 60H30, 93E05, 90C39, 35K55, 35K65}

\vspace{5mm}\noindent \textbf{Keywords.} BSDE, HJB equation, Lipschitz continuity, reflected BSDE, semi-concavity, value function.

\newpage
\section{\large{Introduction}}

\vspace{2mm}

We are interested in regularity properties of possibly degenerate Hamilton-Jacobi-Bellman (HJB) equations with or without obstacles. More precisely, we consider the following HJB equation
\begin{equation}\label{HJB-1}
\frac{\partial}{\partial t}V(t,x)+\inf_{u\in U}H(t,x,V(t,x),\nabla V(t,x),D^2V(t,x),u)=0,
\end{equation}
and the following HJB equation with either a lower obstacle
\begin{equation}\label{HJB-2}
\min\left\{V(t,x)-\varphi(t,x),-\frac{\partial}{\partial t}V(t,x)-\inf_{u\in U}H(t,x,V(t,x),\nabla V(t,x),D^2V(t,x),u)\right\}=0,
\end{equation}
or with an upper obstacle, \begin{equation}\label{HJB-11}
\max\left\{V(t,x)-\varphi(t,x),-\frac{\partial}{\partial t}V(t,x)-\inf_{u\in U}H(t,x,V(t,x),\nabla V(t,x),D^2V(t,x),u)\right\}=0,
\end{equation}
$(t,x)\in[0,T]\times {\mathbb R}^d$ with terminal condition $V(T,x)=\Phi(x),\, x\in {\mathbb R}^d$, and with the Hamiltonian
$$H(t,x,y,p,A,u)=\frac{1}{2}\mbox{tr}\left(\sigma\sigma^*(t,x,u)A\right)+ b(t,x,u)p+f(t,x,y,p\sigma(t,x,u),u),$$
$(t,x,y,p,A,u)\in[0,T]\times {\mathbb R}^d\times {\mathbb R}\times {\mathbb R}^d\times {\mathbb S}^d\times U$, where ${\mathbb S}^d$\ denotes the space of all symmetric $d\times d$\ matrices, and $U$\ is a compact metric control state space.
If $\sigma\sigma^*(t,x,u)\geq \alpha I_{\mathbb{R}^d}\ (\alpha>0)$, the regularity of the solution of the HJB equation (\ref{HJB-1}) is well studied (see, e.g., Krylov~\cite{K}). Here we are interested in the case of possible degeneracy of $\sigma\sigma^*$.

It is well-known that under continuity and growth assumptions on the coefficients, the HJB equations (\ref{HJB-1}), (\ref{HJB-2}) and (\ref{HJB-11}) have a unique viscosity solution $V\in C_p([0, T]\times \mathbb{R}^d)$, respectively; see, e.g., Buckdahn and Li~\cite{bl},~\cite{bl1}, Wu and Yu~\cite{WY}, Crandall, Ishii, Lions~\cite{cil} (the reader more interested in viscosity solution is referred  to the latter reference). Moreover, if the coefficients $b,\ \sigma,\ f$ are continuous and of linear growth, and if $b(t,.,u),\ \sigma(t,.,u),\ f(t,.,.,.,u)$\ are Lipschitz, uniformly with respect to $t, u$, then
\begin{equation}\label{1.3} \begin{aligned}
& \mbox{(i)} |V(t, x) - V(t, x^{\prime})| \leq C|x -x^{\prime}|, \\
& \mbox{(ii)} |V(t, x) - V(t^{\prime}, x)| \leq C(1 + |x|)\sqrt{|t - t^{\prime}|},\\ \end{aligned}
  \end{equation}
$(t, x),\ (t^{\prime}, x^{\prime})\in [0, T]\times \mathbb{R}^d$, for some constant $C\in \mathbb{R}^+$; see, e.g., Lemma 3.5 and Theorem 3.10 in Buckdahn and Li~\cite{bl}, or Peng~\cite{ps}, for the HJB equations (\ref{HJB-1}); Lemma 3.1 and Theorem 3.2 in Buckdahn and Li~\cite{bl1}, or (ii) from the proof of Proposition 3.12 in Wu and Yu~\cite{WY}, for the HJB equations (\ref{HJB-2}) and (\ref{HJB-11}).
\br Indeed, in~\cite{bl} and~\cite{bl1} stochastic differential games and the viscosity solutions of the associated HJB-Isaacs equations with and without obstacles are studied, but, stochastic control problems and associated HJB equations with and without obstacles can be regarded as a special case, in which the control state space of one of the players is a singleton. Therefore, here we can use the results from~\cite{bl} and~\cite{bl1}.
\er
However, here we are interested in regularity properties of $V(t, x)$\ in $(t, x)$. These regularity properties concern the joint Lipschitz property of $V$\ in $(t, x)$, but also the semiconcavity of $V$\ in $(t, x)$, where the semiconcavity is understood in the following sense (the reader is referred to~\cite{bcq} or~\cite{cs}):
\bde Let $A\subset \mathbb{R}^d$ be an open set and let $f: [0, T]\times A\rightarrow \mathbb{R}^n$. We say
that $f$ is $(C_\delta-)$ semiconcave (with linear modulus) in $A$ if for all $\delta>0$, there exists a constant $C_\delta\geq 0$\ such that, for all $x,\ x'\in A$,\ $t,\ t'\in [0, T-\delta],$\ and for all $\lambda\in [0, 1]$,
\be\lambda f(t, x)+(1-\lambda)f(t', x')\leq f(\lambda(t, x)+(1-\lambda)(t', x'))+C_\delta \lambda(1-\lambda)(|t-t'|^2+|x-x'|^2).\ee
Any constant $C_\delta$ satisfying the above inequality is called a semiconcavity constant for
$f$ in $A$.
\ede

However, one has to be careful here. It turns out, and will be pointed out by counterexamples, that the joint Lipschitz continuity and the semiconcavity don't hold on $[0, T]\times \mathbb{R}^d$, but only on $[0, T-\delta]\times \mathbb{R}^d$, for any $\delta>0$. We emphasize the importance of the semiconcavity of $V$\ on $[0, T-\delta]\times \mathbb{R}^d$, for any $\delta>0$, which has, due to Alexandrov's theorem, the immediate consequence that $V$\ admits a second order expansion with respect to $(t, x)$, in almost every $(t, x)\in [0, T]\times \mathbb{R}^d$. Cannarsa and Sinestrari~\cite{cs} (for $\sigma=0$) showed that these regularity properties are the best ones, which can be expected for Hamilton-Jacobi equations. The Lipschitz continuity and semiconcavity of $V(t, x)$\ in $x$, uniformly with respect to $t\in [0, T]$, have been well-known already for a long time. They are the result of straight-forward computations; see, for instance, Fleming and Soner~\cite{fs}, Peng~\cite{ps}, Ishii and Lions~\cite{il}, Yong and Zhou~\cite{yz}. Buckdahn, Cannarsa and Quincampoix~\cite{bcq} studied recently the joint Lipschitz continuity and semiconcavity of solutions $V(t,x)$ of HJB equation without obstacle when $f(t, x, y, z, u)=f(t, x, u)$\ doesn't depend on $(y, z)$. They used a new technique which is a method of time change in the associated stochastic control problem. In this paper we adapt their method to more general HJB equations and to HJB equations with obstacle by developing an associated approach using backward stochastic differential equations (BSDEs). To be more precise, let $(t, x)\in [0, T]\times \mathbb{R}^d$, and $W=(W_s)_{s\in [t, T]}$\ be a $m$-dimensional Brownian motion with $W_t=0$. By ${\mathbb F}^W=\{{\cal F}_s^W=\sigma\{W_r, r\leq s\}\bigvee {\cal N}_P\}_{s\in [t, T]}$\ we denote the filtration generated by $W$\ and augmented by all $P$-null sets. We consider the following forward stochastic differential equation (SDE)
\begin{equation}
\label{FSDE1.1}\left\{\begin{aligned}
&dX_s^{t,x,u}=\sigma(s, X_s^{t,x,u}, u_s)dW_s + b(s, X_s^{t,x,u},u_s)ds, \quad \quad \quad \quad s\in[t, T],\\
&X_t^{t,x,u}=x,
\end{aligned}\right.\end{equation}
which we associate with the reflected backward stochastic differential equation (RBSDE) with a lower obstacle
\begin{equation}
\label{RBSDE1.2}\left\{\begin{aligned}
&dY_s^{t,x,u}= -f(s, X_s^{t,x,u}, Y_s^{t,x,u}, Z_s^{t,x,u}, u_s)ds + Z_s^{t,x,u}dW_s - dK_s^{t,x,u},\\
&Y_T^{t,x,u}=\Phi(X_T^{t,x,u}),\ K^{t,x,u} \mbox{continuous, increasing},\ K_t^{t,x,u}=0,\\
&Y_s^{t,x,u}\geq \varphi(s, X_s^{t,x,u}), \quad
(Y_s^{t,x,u}-\varphi(s,X_s^{t,x,u}))dK_s^{t,x,u}=0, \quad s\in
[t,T],\end{aligned}\right.\end{equation}
where, the admissible controls $u$\ belong to the space ${\cal U}^W_{t,T}:=L^0_{{\mathbb{F}}^W}(t,T;U)$\ of ${\mathbb{F}}^W$-adapted $U$-valued processes, and $U$ is a compact metric space.
\noindent The coefficients
$$\sigma:[0,T]\times {\mathbb R}^d\times U\rightarrow {\mathbb R}^{d\times m},\, b:[0,T]\times {\mathbb R}^d\times U\rightarrow {\mathbb R}^d,$$
$$ f:[0,T]\times {\mathbb R}^d\times {\mathbb R}\times {\mathbb R}^m\times U\rightarrow {\mathbb R},\, \Phi:{\mathbb R}^d\rightarrow {\mathbb R}\ \mbox{and}\ \varphi:[0,T]\times {\mathbb R}^d\rightarrow {\mathbb R}$$

\noindent are continuous functions which we suppose to satisfy the following standard conditions:

\medskip

H1) The functions $\sigma(.,.,u),b(.,.,u),f(.,.,.,.,u),\varphi(.,.)$ are Lipschitz in $(t,x, y, z)\in[0,T]\times {\mathbb R}^d\times {\mathbb R}\times {\mathbb R}^m,$ uniformly with respect to $u\in U,$ and $\Phi:{\mathbb R}^d\rightarrow {\mathbb R}$ is Lipschitz in $x\in {\mathbb R}^d.$

H2) The functions $\sigma,b,f,\varphi$ and $\Phi$ are bounded.

H3) $\Phi(x)\geq\varphi(T, x), x\in {\mathbb R}^d.$

\medskip
The above RBSDE was introduced in El Karoui, Kapoudjian, Pardoux, Peng and Quenez~\cite{ekppq}. It extends the notion of BSDEs, which was the first time studied in its general form by Pardoux and Peng~\cite{pp}, by endowing it with a lower or an upper barrier.

Then from~\cite{K} and~\cite{ekppq} we know that SDE (\ref{FSDE1.1}) and RBSDE (\ref{RBSDE1.2}) have a unique ${\mathbb F}^W$-adapted, square integrable solution
$X^{t,x,u},$\ and $(Y^{t,x,u}, Z^{t,x,u}, K^{t,x,u})$, respectively. From~\cite{bl1} (or~\cite{WY}) we know that the deterministic function
\begin{equation}\label{Valuefunction1.1}
V(t,x):=\inf_{u\in{\cal U}^W_{t,T}}Y^{t,x,u}_t,\, (t,x)\in[0,T]\times {\mathbb R}^d
\end{equation}
belongs to $C_l([0, T]\times {\mathbb R}^d)$ , and is the unique viscosity solution (unique in $C_p([0, T]\times {\mathbb R}^d)$) of HJB equation (\ref{HJB-2}) with obstacle.
\br By $C_l([0, T]\times {\mathbb R}^d)$\ (respectively, $C_p([0, T]\times {\mathbb R}^d)$) we denote the space of continuous real functions over $[0,T]\times {\mathbb R}^d$\ which have at most linear (respectively, polynomial) growth. \er

 For the proof that $V$\ is deterministic, the reader is referred to Proposition 3.3 in~\cite{bl} or Proposition 3.1 in~\cite{bl1}. Using the time change method in the above control problem for SDE (\ref{FSDE1.1}) and RBSDE (\ref{RBSDE1.2}) we get our main results.

\bt Under the assumptions (H1)-(H3), $V(t,x)$ is joint Lipschitz continuous in $(t, x)\in [0, T-\delta]\times \mathbb{R}^d$, for all $\delta>0,$\  i.e., there exists $C_\delta>0$\ such that, for any $(t, x), (t', x')\in [0, T-\delta]\times \mathbb{R}^d$, \be\label{lip1.1}|V(t, x)-V(t', x')|\leq C_\delta(|t-t'|+|x-x'|).\ee
\et

In fact, we will even show more: the value functions $V_n, n\geq 1,$\ of the associated stochastic control problem in which the reflected BSDE is replaced by the penalized one (see (\ref{valuefunction-n}) and (\ref{lip2.1})), satisfy (\ref{lip1.1}), uniformly with respect to $n\geq 1$.

\br A symmetric argument shows that the continuous viscosity solution $V(t,x)$ of equation (\ref{HJB-11}) with an upper obstacle also satisfies the joint Lipschitz property as that stated in Theorem 1.1 for the viscosity solution of the equation (\ref{HJB-2}) with a lower obstacle. For the stochastic interpretation of the solution $V$ of equation (\ref{HJB-11}) the reader is referred to (\ref{4.2}).
\er

Concerning the joint semiconcavity which is our second main result, we will give a counterexample which shows that the viscosity solution $V$ of HJB equation (\ref{HJB-2}) with a lower obstacle is, in general, not semiconcave on $[0, T-\delta]\times {\mathbb R}^d\ (\delta>0)$, even if the lower obstacle is constant. However, if $V$ is the viscosity solution of HJB equation (\ref{HJB-11}) with an upper obstacle, then $V$ has the joint semiconcavity property in $(t, x)\in [0, T-\delta]\times \mathbb{R}^d$, for all $\delta>0$. For this we need the following assumptions:

H3') $\Phi(x)\leq\varphi(T, x), x\in {\mathbb R}^d.$

H4) $f(t, x, y, z, u)$\ is semiconcave in $(t, x, y, z)\in [0,T]\times {\mathbb R}^d\times {\mathbb R}\times {\mathbb R}^m$, uniformly with respect to $u\in {U}$\ (i.e., the semiconcavity constant $C_\delta$\ doesn't depend on $u$); $\Phi(x)$ is semiconcave.

H5) The first-order derivatives $\nabla_{(t,x)}b,\  \nabla_{(t,x)}\sigma$\ of $b$\ and $\sigma$\ with respect to $(t, x)$ exist, are continuous in $(t, x, u)$ and Lipschitz continuous in $(t,x)$, uniformly with respect to $u\in {U}.$

H6) $f(t, x, y, z, u)=f(t, x, y, u)$\ is independent of $z$; $\varphi$\ is semiconcave in $(t, x)\in [0,T]\times {\mathbb R}^d$.

H7) $\varphi(t, x)=\varphi\in {\mathbb R}$, $(t, x)\in [0,T]\times {\mathbb R}^d$.

\bt In addition to (H1), (H2), and (H3'), we assume that H4), H5) as well as either H6) or H7) hold. Then, the value function $V$\ which is the viscosity solution of HJB equation (\ref{HJB-11}), is $(C_\delta)$-semiconcave on $[0, T-\delta]\times \mathbb{R}^d$, for all $\delta>0$.
\et

\br A standard transformation allows to replace the assumption H7) of constancy of $\varphi$\ by that $\varphi\in C_b^{3,4} ([0,T]\times {\mathbb R}^d)$. For simplicity we restrict ourselves to H7). However, also here for the case of semiconcavity we will prove even more: under the assumptions of the theorem the value functions $V_n,\ n\geq 1,$ of the associated stochastic control problem, in which the reflected BSDE is replaced by penalized ones (see,
(\ref{valuefunction-n}) and (\ref{4.5})), are $C_\delta$-semiconcave on $[0, T-\delta]\times \mathbb{R}^d$, uniformly with respect to $n\geq 1$, for all $\delta>0$.\er

\br 1) The boundedness assumption on the coefficients is made to simplify the computations and to emphasize the main arguments.

2) The above two theorems remain valid for HJB equations (\ref{HJB-1}) without obstacle. Indeed, all coefficients are bounded, and so the viscosity solution $V(t, x)$\ of the HJB equation without obstacle is $|V(t,x)|\leq C,\ (t, x)\in [0, T]\times \mathbb{R}^d$, \ for some $C\in \mathbb{R}$\ depending only on the bounds of $\sigma,\ b,\ f$\ and $\Phi$. It suffices to suppose that the obstacle $\varphi$\ is sufficiently large, i.e., $|\varphi(t, x)|\geq C,\ (t, x)\in [0, T]\times \mathbb{R}^d$, in order to interpret $V$\ as a solution of HJB equation with obstacle. On the other hand, the associated BSDE becomes a RBSDE with a lower obstacle or an upper one, see Remark \ref{remark2.1}. Therefore, we only need to study HJB (\ref{HJB-2}) or (\ref{HJB-11}).
\er

Our paper is organized as follows. In Section 2, we study the joint Lipschitz continuity for the HJB equations with or without obstacles with the help of the associated stochastic control problems which cost functionals are given by BSDEs or by RBSDEs. For this end, a special BSDE method based on a time change is developed. Section 3 is devoted to study the semi-concavity for the HJB equations with or without obstacles. We prove that, under some appropriate assumptions, the viscosity solution $V$ also satisfies the semiconcavity property, but only if it is the solution of a HJB equation (\ref{HJB-11}) with an upper obstacle. Our analysis is based on the combination of two time changes and the development of appropriate BSDE estimates under time change. Concerning the viscosity solution of a HJB equation (\ref{HJB-2}) with a lower obstacle, we show with a simple counter-example that semiconcavity is, in general, not satisfied. For the purpose of readability some basics on BSDEs and RBSDEs are given, but postponed to the Appendix (Section 4).

\section{The joint Lipschitz continuity of the value function}

\vspace{2mm}

Given a compact metric control state space $U$ we consider the Hamilton-Jacobi-Bellman (HJB) equation with a lower obstacle
\begin{equation}\label{HJB-obst}
\min\left\{V(t,x)-\varphi(t,x),-\frac{\partial}{\partial t}V(t,x)-\inf_{u\in U}H(t,x,V(t,x),\nabla V(t,x),D^2V(t,x),u)\right\}=0,
\end{equation}
$(t,x)\in[0,T]\times {\mathbb R}^d$, with terminal condition $V(T,x)=\Phi(x),\, x\in {\mathbb R}^d$, and with the Hamiltonian
$$H(t,x,y,p,A,u)=\frac{1}{2}\mbox{tr}\left(\sigma\sigma^*(t,x,u)A\right)+ b(t,x,u)p+f(t,x,y,p\sigma(t,x,u),u),$$
$(t,x,y,p,A,u)\in[0,T]\times {\mathbb R}^d\times {\mathbb R}\times {\mathbb R}^d\times S^d\times U$.

\noindent The coefficients
$$\sigma:[0,T]\times {\mathbb R}^d\times U\rightarrow {\mathbb R}^{d\times m},\, b:[0,T]\times {\mathbb R}^d\times U\rightarrow {\mathbb R}^d,$$
$$ f:[0,T]\times {\mathbb R}^d\times {\mathbb R}\times {\mathbb R}^m\times U\rightarrow {\mathbb R},\, \Phi:{\mathbb R}^d\rightarrow {\mathbb R}\ \mbox{and}\ \varphi:[0,T]\times {\mathbb R}^d\rightarrow {\mathbb R}$$

\noindent are continuous functions which we suppose to satisfy H1)-H3).

\medskip

It is by now well-known (see, for instance,~\cite{bl1},~\cite{WY}) that the above HJB equation with the obstacle possesses a continuous viscosity solution $V\in C_{b}([0,T]\times {\mathbb R}^d)$\ (the space of bounded continuous functions over $[0,T]\times {\mathbb R}^d$) which is unique in the class of viscosity solutions with polynomial growth. It can be stochastically interpreted by the following controlled stochastic system.

Let $(t, x)\in[0,T]\times {\mathbb R}^d$. Given a $m$-dimensional Brownian motion $W=(W_s)_{s\in[t,T]}$ with $W_t=0,$ defined on a complete probability space $(\Omega,{\cal F},P)$ endowed with the filtration ${\mathbb{F}}^W=({\cal F}_s^W)_{s\in[t,T]}$ generated by the Brownian motion $W$ and completed by all $P$-null sets. We introduce the following spaces which will be used frequently in what follows:
    $${\cal{S}}_{{\mathbb{F}}^W}^2(t, T;{\mathbb{R}}^{d}):=\{(\psi_s)_{t\leq s\leq T}\ {\mathbb{R}}^{d}\mbox{-valued}\ {\mathbb{F}}^W\mbox{-adapted continuous process}: E[\sup\limits_{t\leq s\leq T}| \psi_{s} |^2]< +\infty \}; $$
   $$L_{{\mathbb{F}}^W}^{2}(t,T;{\mathbb{R}}^{d}):=\{(\psi_s)_{t\leq s\leq T}\ {\mathbb{R}}^{d}\mbox{-valued }\ {\mathbb{F}}^W\mbox{-progressively measurable}:\ E[\int^T_t| \psi_s| ^2ds]<+\infty \};$$

\vskip0.2cm

$A_{{\mathbb{F}}^W}^2(t, T):=\{(\psi_s)_{t\leq s\leq T}\in {\cal{S}}_{{\mathbb{F}}^W}^2(t, T;{\mathbb{R}}):\ \psi\ \mbox{increasing process},\ \psi_t=0\}.$
    \vskip0.2cm

We consider the forward stochastic differential equation (SDE)
\begin{equation}
\label{FSDE}\left\{\begin{aligned}
&dX_s^{t,x,u}=\sigma(s, X_s^{t,x,u}, u_s)dW_s + b(s, X_s^{t,x,u},u_s)ds, \quad \quad \quad \quad s\in[t, T],\\
&X_t^{t,x,u}=x,
\end{aligned}\right.\end{equation}
which we associate with the reflected backward stochastic differential equation (RBSDE) with a lower barrier
\begin{equation}
\label{RBSDE}\left\{\begin{aligned}
&dY_s^{t,x,u}= -f(s, X_s^{t,x,u}, Y_s^{t,x,u}, Z_s^{t,x,u}, u_s)ds + Z_s^{t,x,u}dW_s - dK_s^{t,x,u},\\
&Y_T^{t,x,u}=\Phi(X_T^{t,x,u}),\\
&Y_s^{t,x,u}\geq \varphi(s, X_s^{t,x,u}), \quad
(Y_s^{t,x,u}-\varphi(s,X_s^{t,x,u}))dK_s^{t,x,u}=0, \quad s\in
[t,T].\end{aligned}\right.\end{equation}
The control process $u$ runs the set of admissible controls ${\cal U}^W_{t,T}:=L^0_{{\mathbb{F}}^W}(t,T;U)$, defined as set of all ${\mathbb{F}}^W$-progressively measurable processes over $(\Omega,{\cal F},P)$, taking their values in $U$. Then, from~\cite{K} and~\cite{ekppq} we know SDE (\ref{FSDE}) and RBSDE (\ref{RBSDE}) have a unique solution
$$(X^{t,x,u}, (Y^{t,x,u}, Z^{t,x,u}, K^{t,x,u}))\in {\cal{S}}_{{\mathbb{F}}^W}^2(t, T;{\mathbb{R}}^{d})\times {\cal{S}}_{{\mathbb{F}}^W}^2(t, T;{\mathbb{R}})\times L_{{\mathbb{F}}^W}^{2}(t,T;{\mathbb{R}}^{m})\times A_{{\mathbb{F}}^W}^2(t, T) .$$
In order to emphasize that we have to deal with the solution of a decoupled forward-backward system driven by the Brownian motion $W$, we also write
$$(X^{t,x,u}(W), (Y^{t,x,u}(W), Z^{t,x,u}(W), K^{t,x,u}(W)))=(X^{t,x,u}, (Y^{t,x,u}, Z^{t,x,u}, K^{t,x,u})).$$
Observe that $Y_t^{t,x,u}$\ is ${\cal F}_t^W$-measurable, and can, hence, be identified with the deterministic real value $E[Y_t^{t,x,u}]$. Moreover, from~\cite{bl1} or~\cite{WY} we know
\begin{equation}\label{Valuefunction}
V(t,x):=\inf_{u\in{\cal U}^W_{t,T}}Y^{t,x,u}_t(W),\, (t,x)\in[0,T]\times {\mathbb R}^d,
\end{equation}
belongs to $C_b([0, T]\times {\mathbb R}^d)$\ and is the unique viscosity solution (unique in $C_p([0, T]\times {\mathbb R}^d)$) of HJB equation (\ref{HJB-obst}) with the obstacle.
Standard SDE and BSDE estimates allow to show (see, e.g.,~\cite{bl1} or~\cite{WY}) that, for all $t,\ t'\in[0,T],\, x,\ x'\in {\mathbb R}^d$,
  \begin{equation}\label{3.5} \begin{aligned}
& \mbox{(i)} |V(t, x)| \leq C,\\
& \mbox{(ii)} |V(t, x) - V(t, x^{\prime})| \leq C|x -x^{\prime}|, \\
& \mbox{(iii)} |V(t, x) - V(t^{\prime}, x)| \leq C(1 + |x|)\sqrt{|t - t^{\prime}|}.\\ \end{aligned}
  \end{equation}
\br\label{remark2.1} The above constant $C$ depends only on the bounds and the Lipschitz constants of the functions $\sigma,b,f,\varphi$ and $\Phi.$ We also observe that, if the coefficients $f$ and $\Phi$ are bounded by $C_0\in {\mathbb R}^{+}$ and $\varphi\le -C_0(1+T)$ on $[0,T]\times {\mathbb R}^d$, then from the comparison theorem-Lemma 4.4 in Section 4

\smallskip

\centerline{$|Y^{t,x,u}_s|\le C_0(1+T)$, and $K^{t,x,u}_s=0,\ s\in[t,T],\ (t,x)\in [0,T]\times {\mathbb R}^d.$}

\noindent This also shows that, by choosing $\overline{\varphi}(t,x)=-C_0(1+T),\, (t,x)\in[0,T]\times {\mathbb R}^d$, a BSDE with the coefficients $f$ and $\Phi$ bounded by $C_0$ can be regarded as a reflected BSDE with a lower barrier $\overline{\varphi}$, which  coefficients satisfy our standard assumptions of boundedness and Lipschitz continuity. Similarly, by choosing $\widehat{\varphi}(t,x)=C_0(1+T),\, (t,x)\in[0,T]\times {\mathbb R}^d$, a BSDE with the coefficients $f$ and $\Phi$ bounded by $C_0$ can be regarded as a reflected BSDE with an upper barrier $\widehat{\varphi}$, which coefficients satisfy our standard assumptions of boundedness and Lipschitz continuity. That means, the value function defined by (\ref{Valuefunction}) is the unique viscosity solution of HJB equation (\ref{HJB-1}) without obstacle. So our studies of the regularity properties of the solutions of HJB equations with obstacles include in particular those without obstacles.\er

 Unlike (\ref{3.5}) our objective here is to study the joint Lipschitz continuity of $V(t, x)$\ in $(t, x)$. This joint Lipschitz property of the solution $V$\ of such HJB equations was somewhat expected, see, Krylov~\cite{K}. However, it doesn't hold on $[0,T]\times {\mathbb R}^d$\ as the following example shows.

\begin{example} We let the dimension $m=d=1$, and we choose the coefficients $b=0,\ \sigma=1,\ f=0,$\ and $\Phi(x)=|x|,\ x\in {\mathbb R}^d$. Then
$$V(t, 0)=E[\Phi(X_T^{t,0})]=E[|W_T-W_t|]=\sqrt{\frac{2}{\pi}}\sqrt{T-t},\ t\in [0, T].$$
It's obvious that $V$ is not Lipschitz in $t$\ and, hence, not jointly Lipschitz in $(t, x)$\ for $t$\ around $t=T$; however, $V$\ is jointly Lipschitz on $[0,T-\delta]\times {\mathbb R}$, for all $\delta>0$.
\end{example}

Our objective in this section is to investigate the joint Lipschitz property of the value function $V$. More precisely, we have the following
\begin{theorem}\label{Th3.1} Under our standard assumptions H1)-H3) the value function $V(.,.)$ is jointly Lipschitz continuous on $[0, T-\delta]\times \mathbb{R}^{d},$  for all $\delta>0.$\end{theorem}
The proof of this theorem will be split into a sequel of different statements. It formalizes and generalizes the method of time change for the underlying Brownian motion introduced into the frame of stochastic control problems with classical cost functional in~\cite{bcq}.

\smallskip

In order to estimate the reflected BSDEs (\ref{RBSDE}) driven by $W$, we approximate them by penalized BSDEs. More precisely, we approximate (\ref{RBSDE}) with its unique solution $(Y^{t,x,u}, Z^{t,x,u}, K^{t,x,u})$
by the following penalized BSDEs:
\begin{equation}\label{penBSDE-0}\left\{\begin{aligned}&dY_s^{t,x,u;n}= -[f(s, X_s^{t,x,u},Y_s^{t,x,u;n}, Z_s^{t,x,u;n}, u_s)+n\left(Y_s^{t,x,u;n} - \varphi(s, X_s^{t,x,u})\right)^{-}]ds\\
&\ \hskip1.5cm  + Z_s^{t,x,u;n}dW_s,\, s\in[t,T],\\
&Y_T^{t,x,u;n}=\Phi(X_T^{t,x,u}),\ n\geq 1.\end{aligned}\right.\end{equation}

For all $n\geq 1$, BSDE (\ref{penBSDE-0}) has a unique solution $(Y^{t,x,u;n},Z^{t,x,u;n})\in{\cal S}^2_{{\mathbb{F}}^{W}}(t,T)\times L^2_{{\mathbb{F}}^{W}}(t,T;{\mathbb R}^m)$.  We define
\be\label{valuefunction-n}V_n(t,x):=\inf_{u\in{\cal U}^W_{t,T}}Y^{t,x,u;n}_t,\, (t,x)\in[0,T]\times {\mathbb R}^d.
\ee
\begin{proposition}\label{p3.1} Under our standard assumptions H1) and H3) we have\\
\begin{equation}\label{3.8}\begin{aligned}
&{\rm{i)}}\ \ Y^{t,x,u;1}_s\le Y^{t,x,u;2}_s\le\cdots \le Y^{t,x,u;n}_s\rightarrow Y^{t,x,u}_s, \mbox{as}\ n\rightarrow+\infty,\, s\in[t,T],\ \mbox{P-a.s.},\ u\in{\cal U}^W_{t,T};\\
&{\rm{ii)}}\ \ \displaystyle E\hskip -0.1cm\left[\sup_{s\in[t,T]}\hskip -0.1cm\left|Y^{t,x,u;n}_s-Y_s^{t,x,u}\right|^2\hskip -0.1cm+\hskip -0.1cm\int_{t}^T\hskip -0.1cm\left|Z^{t,x,u;n}_s-Z^{t,x,u}_s\right|^2ds\right.
\\
&\mbox{ }\ \ \ \  \ \left.\hskip -0.1cm+\hskip -0.1cm\sup_{s\in[t,T]}\hskip -0.1cm\left|K^{t,x,u}_s-n\int_t^s(Y^{t,x,u;n}_r-\varphi(r,X_r^{t,x,u}))^-dr\right|^2\right]\rightarrow 0, \mbox{as}\  n\rightarrow +\infty,\ u\in{\cal U}^W_{t,T};\\
&{\rm{iii)}}\ V_1(t,x)\leq V_2(t, x)\le\cdots \le V_n(t, x)\rightarrow V(t, x), \mbox{as}\ n\rightarrow+\infty,\  (t, x)\in [0,T]\times {\mathbb R}^d.
\end{aligned}\end{equation}
\end{proposition}

\noindent For the proof of these classical results, in particular those of i) and ii), the reader is referred to Section 6 of~\cite{ekppq}. The result iii) can be consulted, for instance, in Theorem 4.2 in~\cite{bl1} or in Lemma 4.3 in~\cite{WY}.
\smallskip

Theorem \ref{Th3.1} follows from the following theorem combined with Proposition \ref{p3.1} iii).
\bt\label{Th3.2} Under the assumptions H1)-H3), for all $\delta>0$, $V_n(t,x)$\ is jointly Lipschitz continuous in $(t, x)\in [0, T-\delta]\times \mathbb{R}^d$, uniformly with respect to $n\geq 1$, i.e., for all $\delta>0,$\ there exists $C_\delta>0$\ such that, for any $n\geq 1,$ $ (t_0, x_0),\ (t_1, x_1)\in [0, T-\delta]\times \mathbb{R}^d$, \be\label{lip2.1}|V_n(t_0, x_0)-V_n(t_1, x_1)|\leq C_\delta(|t_0-t_1|+|x_0-x_1|).\ee
\et

The proof is based on the method of time change and split into several steps. Let us arbitrarily fix $\delta>0$, $ (t_0, x_0),\ (t_1, x_1)\in [0, T-\delta]\times \mathbb{R}^d$. Moreover, let $W^0=(W_s^0)_{s\in [t_0, T]}$\ be a m-dimensional Brownian motion with $W^0_{t_0}=0$, and let $u^0\in {\cal U}^{W^0}_{t_0,T}$. With the notations introduced above we put:
\be\label{3.10}(X^0,Y^0,Z^0,K^0)=(X^{t_0,x_0,u^0}(W^0), Y^{t_0,x_0,u^0}(W^0), Z^{t_0,x_0,u^0}(W^0), K^{t_0,x_0,u^0}(W^0)),\ee
\noindent (the unique solution of SDE (\ref{FSDE}) and RBSDE (\ref{RBSDE}) driven by the Brownian motion $W^0$\ and with initial data ($t_0, x_0$)), and
\be\label{3.11}(Y^{0, n}, Z^{0, n})=(Y^{t_0,x_0,u^0; n}(W^0), Z^{t_0,x_0,u^0;n}(W^0)),\ee
(the unique solution of BSDE (\ref{penBSDE-0}) driven by the Brownian motion $W^0$\ and with ($t_0, x_0, u^0$)\ instead of ($t, x, u$)).

We introduce the linear time change $\tau:[t_1,T]\rightarrow [t_0,T]$\ by setting
\be\label{3.12}\tau(s) = t_0 + \frac{T - t_0}{T - t_1}(s - t_1),\ \  s\in[t_1,T],\ee

\noindent and we remark that $\displaystyle\dot{\tau} \left(=\frac{d}{ds}\tau(s)\right)=\frac{T - t_0}{T - t_1}.$\
Consequently,
\be\label{3.13} W_s^1:= W_{\tau(s)}^0\frac{1}{\sqrt{\dot{\tau}}},\,   s\in[t_1,T],\ee
defines a (m-dimensional) Brownian motion with $W_{t_1}^1=0.$ Then, obviously, the time transformed control process $u^1_s:=u^0_{\tau(s)},\, s\in[t_1, T],$ is an admissible control process with respect to the natural filtration ${\mathbb{F}}^{W^1}=({\cal F}_s^{W^1})_{s\in[t_1, T]}$\ generated by the Brownian motion $W^1$ over the time interval $[t_1,T]$: $u^1=(u_s^1)_{s\in[t_1, T]}\in{\cal U}^{W^1}_{t_1, T}\left(=L^0_{{\mathbb{F}}^{W^1}}(t_1,T;U)\right).$

Having a Brownian motion $W^1=(W^1_s)_{s\in[t_1,T]}$ over the time interval $[t_1,T]$ and an associated admissible control $u^1\in{\cal U}^{W^1}_{t_1,T}$ we can solve the corresponding system (\ref{FSDE})-(\ref{RBSDE}), but now driven by the Brownian motion $W^1,$ with $((t_1,x_1),W^1,u^1)$ at the place of $((t_0,x_0),W^0,u^0)$, and we denote its unique solution by
\be\label{3.14}(X^1,Y^1,Z^1,K^1)=(X^{t_1,x_1,u^1}(W^1),Y^{t_1,x_1,u^1}(W^1),Z^{t_1,x_1,u^1}(W^1),K^{t_1,x_1,u^1}(W^1)).\ee
Correspondingly, the solution of the penalized BSDE (\ref{penBSDE-0}) driven by the Brownian motion $W^1$ is denoted by
\be\label{3.15}(Y^{1,n}, Z^{1,n})=(Y^{t_1,x_1,u^1;n}(W^1),Z^{t_1,x_1,u^1;n}(W^1)),\ee
while the associated solution of the forward equation is again $X^1=X^{t_1,x_1,u^1}(W^1).$

Therefore, the above procedure has provided two different forward equations, that for $X^0=X^{t_0,x_0,u^0}(W^0)$ and $X^1=X^{t_1,x_1,u^1}(W^1),$\ i.e.,
\begin{equation}\label{FSDE-0}\begin{aligned}
&dX_s^0=\sigma(s, X_s^0, u^0_s)dW^0_s + b(s, X_s^0,u^0_s)ds,\, s\in[t_0, T],\, X_{t_0}^0=x_0,
\end{aligned}\end{equation}
\begin{equation}\label{FSDE-1}\begin{aligned}
&dX_s^1=\sigma(s, X_s^1, u_s^1)dW^1_s + b(s, X_s^1,u_s^1)ds,\, s\in[t_1, T],\, X_{t_1}^1=x_1,
\end{aligned}\end{equation}
which we associate with the respective RBSDEs
\begin{equation}
\label{RBSDE-0}\left\{\begin{aligned}
&dY_s^{0}= -f(s, X_s^{0}, Y_s^{0}, Z_s^{0}, u^0_s)ds + Z_s^{0}dW^0_s - dK_s^{0},\\
&Y_T^{0}=\Phi(X_T^{0}),\\
&Y_s^{0}\geq \varphi(s, X_s^{0}), \quad
(Y_s^{0}-\varphi(s,X_s^{0}))dK_s^{0}=0, \quad s\in
[t_0,T],\end{aligned}\right.\end{equation}
and
\begin{equation}
\label{RBSDE-1}\left\{\begin{aligned}
&dY_s^{1}= -f(s, X_s^{1}, Y_s^{1}, Z_s^{1}, u^{1}_s)ds + Z^{1}_sdW^1_s - dK_s^{1},\\
&Y_T^{1}=\Phi(X_T^{1}),\\
&Y_s^{1}\geq \varphi(s, X_s^{1}), \quad
(Y_s^{1}-\varphi(s,X_s^{1}))dK_s^{1}=0, \quad s\in
[t_1,T].\end{aligned}\right.\end{equation}
On the other hand, RBSDE (\ref{RBSDE-0}) with its unique solution $(Y^0,Z^0,K^0)$\ is approximated
by the following penalized BSDEs:
\begin{equation}\label{penBSDE-1}\left\{\begin{aligned}&dY_s^{0,n}= -[f(s, X_s^0,Y_s^{0,n}, Z_s^{0,n}, u_s^0)+ n\left(Y_s^{0,n} - \varphi(s, X_s^0)\right)^{-}]ds +
Z_s^{0,n}dW_s^0,\\
&Y_T^{0,n}=\Phi(X_T^0),\ s\in[t_0,T],\ n\geq 1.\end{aligned}\right.\end{equation}
And RBSDE (\ref{RBSDE-1}) with its unique solution $(Y^1, Z^1, K^1)$\ is approximated by the following penalized equations
\begin{equation}\label{penBSDE-2}\left\{\begin{aligned}&dY_s^{1,n}= -[f(s, X_s^1,
Y_s^{1,n}, Z_s^{1,n}, u_s^1)+ n\left(Y_s^{1,n} - \varphi(s, X_s^1)\right)^{-}]ds +Z_s^{1,n}dW_s^1,\\
&Y_T^{1,n}=\Phi(X_T^1),\ s\in[t_1, T],\ n\geq 1.\end{aligned}\right.\end{equation}

In order to be able to compare the both SDEs (\ref{FSDE-0}) and (\ref{FSDE-1}) which are defined over different time intervals and driven by different Brownian motions, we have to make the inverse time change $\tau^{-1}:[t_0,T]\rightarrow[t_1,T],\ \tau^{-1}(s)=t_1+\frac{T-t_1}{T-t_0}(s-t_0),\ s\in [t_0, T],$ in equation (\ref{FSDE-1}) in order to have two SDEs driven by the same Brownian motion $W^0=(W^0_s)_{s\in[t_0,T]}.$ For this we define
\be\label{3.22}\widetilde{{X}}_s^1:= {X}_{\tau^{-1}(s)}^1,\  {\widetilde{Y}}_s^{1,n}:={Y}_{\tau^{-1}(s)}^{1,n},\ \widetilde{{Z}}_s^{1,n}:=\frac{1}{\sqrt{\dot{\tau}}}{Z}_{\tau^{-1}(s)}^{1,n},\  s\in[t_0,T].\ee

By observing that
\be\label{3.23}W_{\tau^{-1}(s)}^1 =\frac{1}{\sqrt{\dot{\tau}}}W^0_s, \mbox{ and } u_{\tau^{-1}(s)}^{1}= u^0_s,\ s\in[t_0,T],\ee
we deduce from (\ref {FSDE-1}) that $\widetilde{{X}}^1=(\widetilde{{X}}_s^1)_{s\in[t_0,T]}$ is the unique continuous ${\mathbb{F}}^{W^0}$-adapted solution of the SDE
\begin{equation}\label{FSDE-2}\begin{aligned}
d\widetilde{X}_s^1=\frac{1}{\sqrt{\dot{\tau}}}\sigma(\tau^{-1}(s), \widetilde{X}_s^1,u^0_s)dW^0_s + \frac{1}{\dot{\tau}} b(\tau^{-1}(s),\widetilde{X}_s^1,u^0_s)ds,\, s\in[t_0, T],\ \ \widetilde{X}_{t_0}^1 = x_1,
\end{aligned}\end{equation}
and from (\ref{penBSDE-2}) we get that $(\widetilde{Y}^{1,n},\widetilde{Z}^{1,n})=(\widetilde{Y}^{1,n}_s,\widetilde{Z}^{1,n}_s)_{s\in[t_0,T]}$ is the unique solution of the penalized BSDE
\begin{equation}\label{penBSDE-3}\left\{\begin{aligned}
&d\widetilde{Y}_s^{1,n} = -\frac{1}{\dot{\tau}}\left[f(\tau^{-1}(s),
\widetilde{X}_s^1, \widetilde{Y}_s^{1,n}, \sqrt{\dot{\tau}}\widetilde{Z}_s^{1,n},
u^0_s) + n(\widetilde{Y}_s^{1,n} - \varphi(\tau^{-1}(s),
\widetilde{X}_s^1))^{-}\right]ds \\
&\mbox{ }\hskip1.1cm+\widetilde{Z}_s^{1,n}dW^0_s,\ \ s\in[t_0,T],\\
&\widetilde{Y}_T^{1,n}= \Phi(\widetilde{X}_T^1),\ n\geq 1.\hfill \end{aligned}\right.\end{equation}

We will prove the following crucial result:
\bp\label{p3.2} There is some $C_\delta\in {\mathbb R}$\ only depending on $\delta$, and on the bounds and the Lipschitz constants of the coefficients such that, for all $n\geq 1,\ s\in [t_0, T],$\ P-a.s.,

\noindent{\rm i)} \be\label{3.26} |\widetilde{Y}_s^{1,n}-Y_s^{0,n}|\leq C_\delta(|t_0-t_1|+\sup_{r\in[t_0, s]}|X_r^0-\widetilde{X}_r^1|). \ee
In particular,
 \be\label{3.27} |{Y}_{t_1}^{1,n}-Y_{t_0}^{0,n}|=|\widetilde{Y}_{t_0}^{1,n}-Y_{t_0}^{0,n}|\leq C_\delta(|t_0-t_1|+|x_0-x_1|). \ee
{\rm ii)} If, in addition, $\varphi(t, x)=\varphi\in {\mathbb R}$, $(t, x)\in [0,T]\times {\mathbb R}^d$, then for all $p\geq 1$\ there is some constant $C_{\delta, p}$\ such that, for all $n\geq 1$, $s\in [t_0, T]$, P-a.s.,
\be\label{3.28}\displaystyle  E[\left(\int_s^T|\widetilde{Z}_r^{1,n}-Z_r^{0,n}|^2dr\right)^p|{\cal F}_s^{W^0}]\leq C_{\delta, p}\left(|t_0-t_1|^2+\sup_{r\in[t_0, s]}|X_r^0-\widetilde{X}_r^1|^2\right)^p. \ee
\ep
Let us begin by showing that Proposition \ref{p3.2} allows to prove Theorem \ref{Th3.2}.

\noindent \textbf{Proof (of Theorem \ref{Th3.2}).} Let $n\geq 1$, and recall that $$\label{valuefunction-n-1}V_n(t_0,x_0):=\inf_{u^0\in{\cal U}^{W^0}_{t_0,T}}Y^{t_0,x_0,u^0;n}_{t_0}=\inf_{u^0\in{\cal U}^{W^0}_{t_0,T}}Y^{0,n}_{t_0}.
$$
Thus, choosing $\epsilon>0$\ arbitrarily small we can find some control $u^0 \in {\cal U}^{W^0}_{t_0,T}$\ (depending on $\epsilon>0$\ and on $n\geq 1$) such that
$$Y^{0,n}_{t_0}\leq V_n(t_0,x_0)+\epsilon.$$
On the other hand, $$\widetilde{Y}_{t_0}^{1,n}=Y_{t_1}^{1, n}\geq V_n(t_1, x_1).$$
Hence, from Proposition \ref{p3.2} we get
\be\label{3.29}\begin{array}{lll}
&V_n(t_1, x_1)-V_n(t_0,x_0)\leq \widetilde{Y}_{t_0}^{1,n}-Y^{0,n}_{t_0}+\epsilon\leq |\widetilde{Y}_{t_0}^{1,n}-Y^{0,n}_{t_0}|+\epsilon\\
&\leq C_\delta(|t_0-t_1|+|X_{t_0}^0-\widetilde{X}_{t_0}^1|)+\epsilon=C_\delta(|t_0-t_1|+|x_0-x_1|)+\epsilon.
\end{array}\ee
Then, the arbitrariness of $\epsilon>0$\ yields $V_n(t_1, x_1)-V_n(t_0,x_0)\leq C_\delta(|t_0-t_1|+|x_0-x_1|),$\ and from the symmetry of the argument we obtain$$|V_n(t_1, x_1)-V_n(t_0,x_0)|\leq C_\delta(|t_1-t_0|+|x_1-x_0|).$$
Finally, by recalling that the constant $C_\delta$\ from Proposition \ref{p3.2} is independent of $(t_0, x_0),\ (t_1, x_1),$\ and $n\geq 1,$\ we complete the proof.\endpf

The proof of Proposition \ref{p3.2} is split into a sequel of lemmas. The first one concerns the comparison of the SDEs (\ref{FSDE-0}) and (\ref{FSDE-2}), i.e.,
$$ \begin{aligned}
&dX_s^0=\sigma(s, X_s^0, u^0_s)dW^0_s + b(s, X_s^0,u^0_s)ds,\ s\in[t_0, T],\ X_{t_0}^0=x_0,
\end{aligned}$$
and,$$ \begin{aligned}
d\widetilde{X}_s^1=\frac{1}{\sqrt{\dot{\tau}}}\sigma(\tau^{-1}(s), \widetilde{X}_s^1,u^0_s)dW^0_s + \frac{1}{\dot{\tau}} b(\tau^{-1}(s),\widetilde{X}_s^1,u^0_s)ds,\ s\in[t_0, T],\ \ \widetilde{X}_{t_0}^1 = x_1.
\end{aligned}$$
To estimate the difference of solutions of these both SDEs, the following lemma turns out to be useful. It can be got by a straight--forward computation (see also~\cite{bcq}).

\begin{lemma}\label{lemma3.2} For the above introduced time change $\tau:[t_1,T]\rightarrow[t_0,T]$ we have
\be\label{3.32}|\tau^{-1}(s)-s|+ \left|\frac{1}{\dot{\tau}}-1\right|+\left|\frac{1}{\sqrt{\dot{\tau}}}-1\right| \le C_\delta|t_0-t_1|,\ s\in [t_0,T],\ee
where the constant $C_\delta$ only depends on $T$ and $\delta>0$, but not on $t_0,\ t_1\in[0, T-\delta].$
\end{lemma}
The above lemma combined with SDE standard estimates allows to get the following result.
\bl
\label{lemma3.1}There is some $C_{\delta, p}\in \mathbb{R}^+$\ only depending on the bounds of $\sigma,\ b$, their Lipschitz constants, and on $\delta,\ p\geq 1$, such that, for all $s\in [t_0, T]$,
\be\label{3.30}
E[\sup_{r\in[s, T]}|X_r^0-\widetilde{X}_r^1|^p|{\cal F}_s^{W^0}]\leq C_{\delta,p}(|t_0-t_1|^p+|X_s^0-\widetilde{X}_s^1|^p),\ \ \mbox{P-a.s.}\ee
In particular, for $s=t_0$,\be\label{3.31}
E[\sup_{r\in[t_0, T]}|X_r^0-\widetilde{X}_r^1|^p]\leq C_{\delta, p}(|t_0-t_1|^p+|x_0-x_1|^p).\ee
\el

\noindent{\textbf{Proof.}} Taking the difference between the SDEs (\ref{FSDE-0}) and (\ref{FSDE-2}) we obtain
\begin{equation}\label{3.33}\begin{aligned}
d(X_s^0-\widetilde{X}_s^1)=&\left(\sigma(s, X_s^0, u^0_s)-\frac{1}{\sqrt{\dot{\tau}}}\sigma(\tau^{-1}(s), \widetilde{X}_s^1,u^0_s)\right)dW^0_s\\
&+ \left(b(s, X_s^0,u^0_s)-\frac{1}{{\dot{\tau}}}b(\tau^{-1}(s),\widetilde{X}_s^1,u^0_s)\right)ds,\, s\in[t_0, T],\\ X_{t_0}^0-\widetilde{X}_{t_0}^1=&x_0-x_1.
\end{aligned}\end{equation}
Thus, taking into account that $b$\ and $\sigma$\ are bounded, SDE standard estimates yield that, for all $p\ge 1$ there is some constant $C_p$ only depending on the bounds and the Lipschitz coefficients of $\sigma$ and $b$ as well as of $T$ and $p$, such that
\be\label{3.34}\mathbb{E}[\sup_{s\leq r\leq T}|\widetilde{X}_r^1 - X_r^0|^p|\mathcal{F}_s^{W^0}] \leq C_{p} \cdot \left(\left|\frac{1}{\dot{\tau}}-1\right|^p+\left|\frac{1}{\sqrt{\dot{\tau}}}-1\right|^p +\int_s^T|\tau^{-1}(r)-r|^pdr+|\widetilde{X}_s^1 - X_s^0|^p\right),\ee
$t_0\leq s\leq T,\, p\geq 1.$ Finally, by applying the preceding lemma we complete the proof.

\bigskip

We will also need the following lemma.
\bl\label{lemma3.3}
${\rm{i)}}\ \mbox{There exists some constant}\ C\ \mbox{only depending on the bounds of}\ f,\ \Phi\ \mbox{and}\ \varphi,$ $ \mbox{such that}$  \begin{equation}\label{3.35}
\left|{Y}_s^{i,n}\right|\le C,\ s\in[t_i,T],\  n\ge 1,\ i=0, 1,\ \mbox{P-a.s.}\end{equation}
{\rm ii)} For all $p\geq 1$\ there is some constant $C_p$\ only depending on the bounds of the coefficients $f,\ \Phi\ \mbox{and}\ \varphi,$  and on $p$, such that $s\in[t_i,T],\  n\ge 1,\ \ i=0, 1,$
\begin{equation}\label{3.36}\begin{aligned}\displaystyle E\left[\left(\int_{s}^T\left|{Z}^{i,n}_r\right|^2dr\right)^p+
\left(n\int_s^T({Y}^{i,n}_r-\varphi(r,X_r^i))^-dr\right)^{2p} \|{\cal F}^{W^i}_s\right]\le C_p,\   \mbox{P-a.s.}\end{aligned}\end{equation}
\el

\noindent\textbf{Proof.} Assertion i) follows directly from Proposition 2.1-i), and the comparison theorem for reflected BSDEs (Lemma \ref{comparisonrbsde} in Section 4) and the boundedness of the coefficients $f,\ \Phi,\ \mbox{and}\ \varphi$.

ii) From the penalized BSDEs (\ref{penBSDE-1}) and (\ref{penBSDE-2}), i) and the boundedness of the coefficients $f$\ and $\Phi$\ we have, for some constant $C_p$,
$$n\int_s^T({Y}^{i,n}_r-\varphi(r,X_r^i))^-dr\leq C_p+\int_s^TZ_r^{i, n}dW_r^i,\ s\in [t_i, T],\ n\geq 1.$$
Hence, \be\label{3.37}E[(n\int_s^T({Y}^{i,n}_r-\varphi(r,X_r^i))^-dr)^{2p}|{\cal F}_s^{W^i}]\leq C_p+C_pE[(\int_s^T|Z_r^{i, n}|^2dr)^p|{\cal F}_s^{W^i}],\ s\in[t_i, T].\ee

On the other hand, from It\^{o}'s formula:
\be\label{3.38}\begin{aligned}\displaystyle
&|Y_s^{i, n}|^2+\int_s^T|Z_r^{i, n}|^2dr=|\Phi(X_T^i)|^2+2\int_s^TY_r^{i,n}f(r, X_r^i, Y_r^{i,n}, Z_r^{i,n}, u_r^i)dr\\
&+2n\int_s^TY_r^{i,n}({Y}^{i,n}_r-\varphi(r,X_r^i))^-dr-2\int_s^TY_r^{i,n}Z_r^{i, n}dW_r^i,\ s\in[t_i, T],\ n\geq 1.\\
\end{aligned}
\ee
From (\ref{3.38}) together with i),
\be\label{3.39}\begin{aligned}\displaystyle
E[(\int_s^T|Z_r^{i, n}|^2dr)^{2p}|{\cal F}_s^{W^i}]\leq & C_p+C_pE[(n\int_s^T({Y}^{i,n}_r-\varphi(r,X_r^i))^-dr)^{2p}|{\cal F}_s^{W^i}]\\
&+C_pE[(\int_s^T|Z_r^{i, n}|^2dr)^p|{\cal F}_s^{W^i}],\ s\in[t_i, T],\ n\geq 1.\\
\end{aligned}
\ee
The result ii) for $p=1$ (see page 719-720 in Section 6 in~\cite{ekppq}) combined with (\ref{3.38}) and (\ref{3.39}) yields the general result ii).\endpf

For the proof of Proposition \ref{p3.2} we have to compare the BSDEs (\ref{penBSDE-1}) and (\ref{penBSDE-3}), i.e.,
\begin{equation}\label{3.40}\left\{\begin{aligned}&dY_s^{0,n}= -[f(s, X_s^0,Y_s^{0,n}, Z_s^{0,n}, u_s^0)+ n\left(Y_s^{0,n} - \varphi(s, X_s^0)\right)^{-}]ds +
Z_s^{0,n}dW_s^0,\\
&Y_T^{0,n}=\Phi(X_T^0),\ \ \ \ \  s\in[t_0,T],\end{aligned}\right.\end{equation}and
\begin{equation}\label{3.40-1}\left\{\begin{aligned}
&d\widetilde{Y}_s^{1,n} = -\frac{1}{\dot{\tau}}\left[f(\tau^{-1}(s),
\widetilde{X}_s^1, \widetilde{Y}_s^{1,n}, \sqrt{\dot{\tau}}\widetilde{Z}_s^{1,n},
u^0_s) + n(\widetilde{Y}_s^{1,n} - \varphi(\tau^{-1}(s),
\widetilde{X}_s^1))^{-}\right]ds \\
&\mbox{ }\hskip1.1cm+\widetilde{Z}_s^{1,n}dW^0_s,\ \ s\in[t_0,T],\\
&\widetilde{Y}_T^{1,n}= \Phi(\widetilde{X}_T^1).\hfill \end{aligned}\right.\end{equation}

But, the different structure of the penalization terms and different obstacles don't allow a direct estimate to get Proposition \ref{p3.2}; so intermediate steps are necessary.

Let us first compare BSDE (\ref{3.40-1}) with the following BSDE (\ref{penBSDE-4}):
\begin{equation}\label{penBSDE-4}
\left\{\begin{aligned}
d\widehat{Y}_s^{1,n} =&-\big[f(s, X_s^0,
\widehat{Y}_s^{1,n} - A_s, \widehat{Z}_s^{1,n}, u^0_s) + n(\widehat{Y}_s^{1,n} -
\varphi(s, X_s^0) - A_s)^{-}\\
&+C_{\delta}|t_0 - t_1|(1 +|\widehat{Z}_s^{1,n}|)+ CA_s\big]ds +\widehat{Z}_s^{1,n}dW^0_s,\ s\in [t_0, T],\\
\widehat{Y}_T^{1,n}=&\Phi(X_T^0)+A_T,\end{aligned}\right.\end{equation}
where $C_\delta,\ C\geq 1$\ are constants which are large enough (their precise choice becomes clear from the proof of the lemma below), and
\be\label{3.43} A_s:= C_{\delta}\sup_{r\in[t_0,s]}(|t_0 - t_1| + |\widetilde{X}_r^1 - X_r^0|),\ s\in[t_0,T].\ee
Note that $A=(A_s)_{s\in [t_0, T]}$\ is an ${\mathbb F}^{W^0}$-adapted, continuous increasing process, $A_{t_0}= C_{\delta}(|t_0 - t_1| + |x_0 - x_1|),$\ and from Lemma \ref{lemma3.1} we see that: for all $q\geq 1$,
\be\label{3.44}\mathbb{E}[A_T^q-A_s^q|{\cal F}_s^{W^0}]\leq C_{\delta, q} \cdot \left(\left|t_0-t_1\right|^q +|\widetilde{X}_s^1 - X_s^0|^q\right), \quad s\in [t_0, T].\ee
\bl\label{lemma3.4} Under our standard assumptions H1)-H3) we have
\be\label{3.45}\widetilde{Y}_s^{1,n}\leq \widehat{Y}_s^{1,n},\ s\in[t_0, T],\ n\geq 1,\  \mbox{P-a.s}.\ee
\el
\noindent \textbf{Proof}. The proof is based on the comparison theorem for BSDEs (Lemma \ref{comparisonbsde} in Section 4). For this we note that, since $\varphi$\ is bounded and Lipschitz,
\be\label{3.46}\begin{aligned}\displaystyle
&|\frac{1}{\dot{\tau}}\varphi(\tau^{-1}(s), \widetilde{X}_s^1)-\varphi(s, X_s^0)|\leq C(|1-\frac{1}{\dot{\tau}}|+|\tau^{-1}(s)-s|+|\widetilde{X}_s^1 - X_s^0|)\\
&\leq C_\delta(|t_1-t_0|+|\widetilde{X}_s^1 - X_s^0|),\ s\in[t_0, T],
\end{aligned}
\ee
(recall Lemma \ref{lemma3.2}).
Thus, recalling that $\widetilde{Y}^{1,n}$ s bounded, uniformly w.r.t. $n\ge 1,$ we get from Lemma \ref{lemma3.2} that
\be\label{3.47}\begin{array}{rcl}
&\widetilde{Y}_s^{1,n}-\varphi(s, X_s^0)\leq \frac{1}{\dot{\tau}}(\widetilde{Y}_s^{1,n}-\varphi(\tau^{-1}(s), \widetilde{X}_s^1))+C_\delta(|t_1-t_0|+|\widetilde{X}_s^1 - X_s^0|)\\
&\leq \frac{1}{\dot{\tau}}(\widetilde{Y}_s^{1,n}-\varphi(\tau^{-1}(s), \widetilde{X}_s^1))+A_s,\ s\in[t_0, T].
\end{array}
\ee
Then, from (\ref{3.47}),
\be\label{3.48} (\frac{1}{\dot{\tau}}(\widetilde{Y}_s^{1,n}-\varphi(\tau^{-1}(s), \widetilde{X}_s^1)))^-\leq (\widetilde{Y}_s^{1,n}-\varphi(s, X_s^0)-A_s)^-,\ s\in[t_0, T].\ee
Moreover, from
\be\label{3.49}\begin{aligned}\displaystyle
&|\frac{1}{\dot{\tau}}f(\tau^{-1}(s),\widetilde{X}_s^1, \widetilde{Y}_s^{1,n}, \sqrt{\dot{\tau}}\widetilde{Z}_s^{1,n}, u_s^0)-f(s, X_s^0, \widetilde{Y}_s^{1,n},\widetilde{Z}_s^{1,n},u_s^0)|\\
&\leq C(|\frac{1}{\dot{\tau}}-1|+|\tau^{-1}(s)-s|+|1-\sqrt{\dot{\tau}}||\widetilde{Z}_s^{1,n}|)+C|\widetilde{X}_s^1-X_s^0|\\
&\leq C_\delta|t_0-t_1|(1+|\widetilde{Z}_s^{1,n}|)+C|\widetilde{X}_s^1-X_s^0|\\
&\leq C_\delta|t_0-t_1|(1+|\widetilde{Z}_s^{1,n}|)+A_s,\ \ s\in[t_0, T],
\end{aligned}
\ee
\noindent we have
\be\label{3.50}\begin{aligned}\displaystyle
&\frac{1}{\dot{\tau}}f(\tau^{-1}(s),\widetilde{X}_s^1, \widetilde{Y}_s^{1,n}, \sqrt{\dot{\tau}}\widetilde{Z}_s^{1,n}, u_s^0)\\
&\leq f(s, X_s^0, \widetilde{Y}_s^{1,n},\widetilde{Z}_s^{1,n}, u_s^0)+C_\delta|t_0-t_1|(1+|\widetilde{Z}_s^{1,n}|)+A_s\\
&\leq f(s, X_s^0, \widetilde{Y}_s^{1,n}-A_s,\widetilde{Z}_s^{1,n}, u_s^0)+C_\delta|t_0-t_1|(1+|\widetilde{Z}_s^{1,n}|)+CA_s,\ s\in[t_0, T].
\end{aligned}
\ee
We also observe that, thanks to the Lipschitz property of $\Phi$,
\be\label{3.51} \Phi(\widetilde{X}_T^1)\leq \Phi({X}_T^0)+A_T,\ \mbox{P-a.s.}\ee

The relations (\ref{3.48}), (\ref{3.50}) and (\ref{3.51}) allow to apply the comparison theorem (Lemma \ref{comparisonbsde} in Section 4) to the both BSDEs, and thus to conclude that
$$\widetilde{Y}_s^{1,n}\leq \widehat{Y}_s^{1,n},\ s\in[t_0, T],\ n\geq 1,\  \mbox{P-a.s}.$$\endpf

The statement of the above lemma can be strengthened as follows:
\bl\label{lemma3.5} Under the standard assumptions H1)-H3) the following holds true:
\be\label{3.52}\begin{aligned}\displaystyle
&{\rm i)}\ -C\leq \widetilde{Y}_s^{1,n}\leq \widehat{Y}_s^{1,n}\leq C_\delta+C_\delta A_s,\ s\in[t_0, T],\ n\geq 1,\  \mbox{P-a.s.};\\
&{\rm ii)}\ E[\int_s^T|\widehat{Z}_r^{1,n}|^2dr|{\cal F}_s^{W^0}]\leq C_\delta(1+A_s^2),\ s\in[t_0, T],\ n\geq 1,\  \mbox{P-a.s.}\end{aligned}\ee\el

\noindent \textbf{Proof}. i) Firstly, from Lemma \ref{lemma3.3} we know that $|\widetilde{Y}_s^{1,n}|\leq C,\ s\in[t_0, T],\ n\geq 1,\  \mbox{P-a.s.}$\  Secondly, thanks to the boundedness of $f$\ and $\varphi$, for some constant $C'$ large enough, we have
\be\label{3.53}\begin{array}{lll}
&f(s, X_s^0,\widehat{Y}_s^{1,n}-A_s, \widehat{Z}_s^{1,n}, u^0_s) + n(\widehat{Y}_s^{1,n}-(\varphi(s, X_s^0)+A_s))^{-} +C_{\delta}|t_0 - t_1|(1 +|\widehat{Z}_s^{1,n}|)+ CA_s\\
&\leq C^{'}+n(\widehat{Y}_s^{1,n}-(C^{'}+C^{'}A_s))^{-} +C_{\delta}(1 +|\widehat{Z}_s^{1,n}|)+ C^{'}A_s,\\
\end{array}
\ee
and $\Phi(X_T^0)+A_T\leq \Phi(X_T^0)+C^{'}A_T$. Hence, we can compare (\ref{penBSDE-4}) with the BSDE (\ref{penBSDE-5}):
\begin{equation}\label{penBSDE-5}
\left\{\begin{aligned}
d\overline{Y}_s^{1,n} =&-\big(C^{'}+n(\overline{Y}_s^{1,n}-(C^{'}+C^{'}A_s))^{-} +C_{\delta}(1 +|\overline{Z}_s^{1,n}|)+ C^{'}A_s\big)ds +\overline{Z}_s^{1,n}dW^0_s,\\
\overline{Y}_T^{1,n}=&\Phi(X_T^0)+C^{'}A_T,\ s\in [t_0, T].\end{aligned}\right.\end{equation}

\noindent From the comparison theorem for BSDEs (Lemma \ref{comparisonbsde} in Section 4) we get that \be\label{3.54}\widehat{Y}_s^{1,n}\leq \overline{Y}_s^{1,n},\  s\in[t_0, T],\ n\geq 1,\  \mbox{P-a.s.}\ee

\noindent On the other hand, putting $\overline{Y}_s^{2,n}:=\overline{Y}_s^{1,n}-C^{'}A_s,\ s\in [t_0, T]$, we get
\begin{equation}\label{penBSDE-6}
\left\{\begin{aligned}
d\overline{Y}_s^{2,n} =&-\big(C^{'}+n(\overline{Y}_s^{2,n}-C^{'})^{-} +C_{\delta}(1 +|\overline{Z}_s^{1,n}|)+ C^{'}A_s\big)ds- C^{'}dA_s+\overline{Z}_s^{1,n}dW^0_s,\\
\overline{Y}_T^{2,n}=&\Phi(X_T^0),\ s\in [t_0, T].\end{aligned}\right.\end{equation}

By observing that,
\be\label{3.64-1}(\overline{Y}_s^{2,n}-C^{'})(\overline{Y}_s^{2,n}-C^{'})^{-}\leq 0,\ee
thanks to It\^{o}'s formula and the boundedness of $\Phi$, for arbitrary $\gamma>0$,
\be\label{3.65-1}\begin{aligned}
&e^{\gamma s}|\overline{Y}_s^{2,n}-C^{'}|^2+E[\int_s^T e^{\gamma r}(\gamma |\overline{Y}_r^{2,n}-C^{'}|^2+|\overline{Z}_r^{1,n}|^2)dr|{\mathcal{F}}_s^{W^0}]\\
& \leq C_{\delta,\gamma}+E[\int_s^T e^{\gamma r}(C_\delta|\overline{Y}_r^{2,n}-C^{'}|^2+\frac{1}{2}|\overline{Z}_r^{1,n}|^2)dr|{\mathcal{F}}_s^{W^0}] \\
&+C_{\gamma}E[A_T^2|{\mathcal{F}}_s^{W^0}]+2E[\int_s^T e^{\gamma r}(\overline{Y}_r^{2,n}-C^{'})C^{'}dA_r|{\mathcal{F}}_s^{W^0}], \ s\in[t_0,T],\ n\geq 1.\\
\end{aligned}\ee
Hence, for $\gamma \geq C_\delta+1$ large enough,
\be\label{3.66-1}\begin{aligned}
&|\overline{Y}_s^{2,n}-C^{'}|^2+E[\int_s^T |\overline{Z}_r^{1,n}|^2dr|{\mathcal{F}}_s^{W^0}]\\
& \leq C_{\delta,\gamma}+C_{\delta,\gamma}A_s^2+\widehat{C}_{\gamma}E[\sup_{r\in[s, T]}|\overline{Y}_r^{2,n}-C^{'}|A_T|{\mathcal{F}}_s^{W^0}], \ s\in[t_0,T],\ n\geq 1,\\
\end{aligned}\ee
where $\widehat{C}_{\gamma}$\ only depends on the coefficients in H1)-H3) and on $\delta,\ \gamma\geq 0$.
Let $1<p<2$ and $q>2$\ be such that $\frac{1}{p}+\frac{1}{q}=1$, and let us choose $\varepsilon>0$ be such that $\widehat{C}_{\gamma}\varepsilon(\frac{2}{2-p})^{\frac{2}{p}}<1$. Then,
\be\label{3.67-1}\begin{array}{lll}
&E[\sup\limits_{r\in[s,T]}|\overline{Y}_r^{2,n}-C^{'}|A_T|\mathcal{F}_s^{W^0}]
\\
& \leq (E[\sup\limits_{r\in[s,T]}|\overline{Y}_r^{2,n}-C^{'}|^p|\mathcal{F}_s^{W^0}])^{\frac{1}{p}}(E[A_T^q|\mathcal{F}_s^{W^0}])^{\frac{1}{q}}\\
& \leq \varepsilon M_{s,t}^{\frac{2}{p}}+\frac{1}{\varepsilon}(E[A_T^q|\mathcal{F}_s^{W^0}])^\frac{2}{q}\\                                      
& \leq \varepsilon M_{s,t}^{\frac{2}{p}}+\frac{1}{\varepsilon}C_{\delta, q}A_s^2,\ \ t_0\leq t\leq s \leq T,\ \ n\geq 1\ \mbox{(see: (\ref{3.43}) and (\ref{3.44}))},
\end{array}\ee
\noindent where
$$M_{s,t}:=E[\sup\limits_{r\in[t,T]}|\overline{Y}_r^{2,n}-C^{'}|^p|{\mathcal{F}}_s^{W^0}].$$
From Doob's martingale inequality, since $\frac{2}{p}>1$,
 \be\label{3.68-1}\begin{array}{lll}
&E[\sup\limits_{s\in[t,T]}M_{s,t}^{\frac{2}{p}}|{\mathcal{F}}_t^{W^0}]\leq (\frac{2}{2-p})^{\frac{2}{p}}E[M_{T,t}^\frac{2}{p}|{\mathcal{F}}_t^{W^0}]\\
&\leq(\frac{2}{2-p})^{\frac{2}{p}}E[\sup\limits_{s\in[t,T]}|\overline{Y}_s^{2,n}-C^{'}|^2|{\mathcal{F}}_t^{W^0}],\ t\in[t_0,T].
\end{array}\ee
Hence, from (\ref{3.66-1}), (\ref{3.67-1}) and (\ref{3.68-1}),
\be\label{3.69-1} \begin{array}{lll}
&E[\sup\limits_{s\in[t,T]}|\overline{Y}_s^{2,n}-C^{'}|^2|{\mathcal{F}}_t^{W^0}]\leq C_{\delta,\gamma}+C_{\delta,\varepsilon}A_t^2\\
& +\widehat{C}_{\gamma}\varepsilon (\frac{2}{2-p})^\frac{2}{p}E[\sup\limits_{s\in[t,T]}|\overline{Y}_s^{2,n}-C^{'}|^2|{\mathcal{F}}_t^{W^0}],\ t\in[t_0,T],
\end{array}\ee
and since $\widehat{C}_{\gamma}\varepsilon (\frac{2}{2-p})^{\frac{2}{p}}<1$, we get
\be\label{3.70-1}E[\sup\limits_{s\in[t,T]}|\overline{Y}_s^{2,n}-C^{'}|^2|{\mathcal{F}}_t^{W^0}]\leq C_{\delta}^{'}(1+A^2_t),\ t\in[t_0,T],\ n\geq 1,\ \mbox{P-a.s.}
\ee
Then we get that $(\overline{Y}_s^{2,n}-C^{'})^2\leq C_\delta^{'}(1+A_s^2),$\ i.e.,
\be\label{3.58} |\overline{Y}_s^{2,n}|\leq C_\delta^{'}(1+A_s),\ s\in [t_0, T]. \ee
Consequently,
$$\widetilde{Y}_s^{1,n}\leq \widehat{Y}_s^{1,n}\leq \overline{Y}_s^{1,n}=\overline{Y}_s^{2,n}+C^{'}A_s\leq C_\delta(1+A_s),\ s\in [t_0, T], \ n\geq 1,\ \mbox{P-a.s.}$$

\noindent ii) Let $C_0$\ be a bound of $\varphi$. Then, from the BSDE (\ref{penBSDE-4}), we have
\begin{equation}\label{3.59}
 \begin{aligned}
&d(\widehat{Y}_s^{1,n}-(C_0+A_s))^2 =\\
&-2(\widehat{Y}_s^{1,n}-(C_0+A_s))\{f(s, X_s^0,\widehat{Y}_s^{1,n}-A_s, \widehat{Z}_s^{1,n}, u^0_s) + n(\widehat{Y}_s^{1,n}-(\varphi(s, X_s^0)+A_s))^{-}\\
 &+C_{\delta}|t_0 - t_1|(1 +|\widehat{Z}_s^{1,n}|)+ CA_s\}ds+|\widehat{Z}_s^{1,n}|^2ds\\
&+2(\widehat{Y}_s^{1,n}-(C_0+A_s))\widehat{Z}_s^{1,n}dW_s^0-2(\widehat{Y}_s^{1,n}-(C_0+A_s))dA_s,\ s\in [t_0, T].\end{aligned} \end{equation}

\noindent Furthermore, from the above result i), $(\widehat{Y}_T^{1,n}-(C_0+A_T))^2\leq C(1+A_T^2)$,
 $(\widehat{Y}^{1,n}_s-(C_0+A_s))(\widehat{Y}_s^{1,n}-(\varphi(s,X_s^0)+A_s))^-\leq0,$
and
 $|\widehat{Y}_s^{1,n}-(C_0+A_s)|\leq C_\delta(1+A_s),\ \ s\in[t_0,T],\ n\geq 1,$
we get by standard estimates
\be\label{3.60}E[\int_s^T|\widehat{Z}_r^{1,n}|^2dr|\mathcal{F}_s^{W^0}]\leq C_\delta+C_\delta E[A_T^2|\mathcal{F}_s^{W^0}]\leq C_\delta+C_\delta(1+A_s^2),\ s\in[t_0,T],\ n\geq 1.\ee
The proof of the lemma is complete.\endpf
Let us put now\be\label{3.61}Y_s^{2,n}:=\widehat{Y}_s^{1,n}-A_s,\ s\in[t_0,T],\ n\geq 1.\ee
Then, from BSDE (\ref{penBSDE-4}) with solution $(\widehat{Y}^{1,n},\widehat{Z}^{1,n})$\ we get
\begin{equation}\label{penBSDE-7}
\left\{
\begin{aligned}
dY_s^{2,n}=&\ -(f(s,X_s^0,Y_s^{2,n},\widehat{Z}_s^{1,n},u_s^0)+n(Y_s^{2,n}-\varphi(s,X_s^0))^{-}\\
           &+C_\delta |t_0-t_1|(1+|\widehat{Z}_s^{1,n}|)+CA_s)ds +\widehat{Z}_s^{1,n}dW_s^0-dA_s,\ s\in[t_0,T],\\
Y_T^{2,n}=&\ \Phi(X_T^0).
\end{aligned}
\right.
\end{equation}
BSDE (\ref{penBSDE-7}) has the advantage that its penalization term is exactly of the same form as that in BSDE (\ref{3.40}). This fact together with the both latter lemmas allow to prove
\bl\label{Lemma3.6} Let us assume H1)-H3). Then, there is some constant $C_\delta$\ such that
\be\label{3.63} E[\sup\limits_{r\in[s,T]}|Y_r^{0,n}-Y_r^{2,n}|^{2}+\int_s^T|Z_r^{0,n}-Z_r^{2,n}|^2dr|\mathcal{F}_s^{W^0}]\leq C_{\delta}A_s^{2},\ee
and, in particular, $|Y_s^{0,n}-Y_s^{2,n}|\leq C_{\delta}A_s$,\  $s\in[t_0,T],\ n\geq 1,$\ P-a.s.
\el
\noindent \textbf{Proof}. We have to compare BSDE (\ref{penBSDE-7}) with BSDE (\ref{3.40}), i.e., with the equation
$$\left\{\begin{aligned}&dY_s^{0,n}= -(f(s, X_s^0,Y_s^{0,n}, Z_s^{0,n}, u_s^0)+n\left(Y_s^{0,n} - \varphi(s, X_s^0)\right)^{-})ds+
Z_s^{0,n}dW_s^0,\, s\in[t_0,T],\\
&Y_T^{0,n}=\Phi(X_T^0).\end{aligned}\right.$$
The proof uses ideas similar to that of (\ref{3.70-1}). However, in view of the importance of the result we prefer to give the proof for the reader's convenience. Taking into account that
\be\label{3.64}(Y_s^{0,n}-Y_s^{2,n})((Y_s^{0,n}-\varphi(s,X_s^0))^--(Y_s^{2,n}-\varphi(s,X_s^0))^-)\leq 0,\ee
we get from standard BSDE estimates that, for arbitrary $\gamma>0$,
\be\label{3.65}\begin{aligned}
&e^{\gamma s}|Y_s^{0,n}-Y_s^{2,n}|^2+E[\int_s^T e^{\gamma r}(\gamma |Y_r^{0,n}-Y_r^{2,n}|^2+|Z_r^{0,n}-\widehat{Z}_r^{1,n}|^2)dr|{\mathcal{F}}_s^{W^0}]\\
& \leq E[\int_s^T e^{\gamma r}(C_\delta|Y_r^{0,n}-Y_r^{2,n}|^2+\frac{1}{2}|Z_r^{0,n}-\widehat{Z}_r^{1,n}|^2)dr|{\mathcal{F}}_s^{W^0}] \\
&+C_{\delta,\gamma}|t_0-t_1|^2E[\int_s^T(1+|\widehat{Z}_r^{1,n}|^2)dr|{\mathcal{F}}_s^{W^0}]
+C_{\delta,\gamma}E[A_T^2|{\mathcal{F}}_s^{W^0}] \\
 &+2E[\int_s^T e^{\gamma r}(Y_r^{0,n}-Y_r^{2,n})dA_r|{\mathcal{F}}_s^{W^0}], \ s\in[t_0,T],\ n\geq 1.\\
\end{aligned}\ee
Hence, for $\gamma \geq C_\delta+1$ large enough,
\be\label{3.66}\begin{aligned}
&|Y_s^{0,n}-Y_s^{2,n}|^2+E[\int_s^T |Z_r^{0,n}-\widehat{Z}_r^{1,n}|^2dr|\mathcal{F}_s^{W^0}]\\
& \leq C_{\delta,\gamma}A^2_s+\widehat{C}_{\delta,\gamma}E[\sup\limits_{r\in[s,T]}|Y_r^{0,n}-Y_r^{2,n}|A_T|{\mathcal{F}}_s^{W^0}],\ s\in[t_0,T],\ n\geq 1,\\
\end{aligned}\ee
where $\widehat{C}_{\delta, \gamma}$\ only depends on the coefficients in H1)-H3) and on $\delta,\ \gamma\geq 0$.
Let $1<p<2$ and $q>2$\ be such that $\frac{1}{p}+\frac{1}{q}=1$, and let us choose $\varepsilon>0$ be such that $\widehat{C}_{\delta, \gamma}\varepsilon(\frac{2}{2-p})^{\frac{2}{p}}<1$. Then,
\be\label{3.67}\begin{array}{lll}
&E[\sup\limits_{r\in[s,T]}|Y_r^{0,n}-Y_r^{2,n}|A_T|\mathcal{F}_s^{W^0}]
\\
& \leq (E[\sup\limits_{r\in[s,T]}|Y_r^{0,n}-Y_r^{2,n}|^p|\mathcal{F}_s^{W^0}])^{\frac{1}{p}}(E[A_T^q|\mathcal{F}_s^{W^0}])^{\frac{1}{q}}\\
& \leq \varepsilon M_{s,t}^{\frac{2}{p}}+\frac{1}{\varepsilon}(E[A_T^q|\mathcal{F}_s^{W^0}])^\frac{2}{q}\\                                      
& \leq \varepsilon M_{s,t}^{\frac{2}{p}}+\frac{1}{\varepsilon}C_{\delta, q}A_s^2,\ \ t_0\leq t\leq s \leq T,\ \ n\geq 1\ \mbox{(see: (\ref{3.43}) and (\ref{3.44}))},
\end{array}\ee
\noindent where
$$M_{s,t}:=E[\sup\limits_{r\in[t,T]}|Y_r^{0,n}-Y_r^{2,n}|^p|{\mathcal{F}}_s^{W^0}].$$
From Doob's martingale inequality, since $\frac{2}{p}>1$,
 \be\label{3.68}\begin{array}{lll}
&E[\sup\limits_{s\in[t,T]}M_{s,t}^{\frac{2}{p}}|{\mathcal{F}}_t^{W^0}]\leq (\frac{2}{2-p})^{\frac{2}{p}}E[M_{T,t}^\frac{2}{p}|{\mathcal{F}}_t^{W^0}]\\
&\leq(\frac{2}{2-p})^{\frac{2}{p}}E[\sup\limits_{s\in[t,T]}|Y_s^{0,n}-Y_s^{2,n}|^2|{\mathcal{F}}_t^{W^0}],\ t\in[t_0,T].
\end{array}\ee
Hence, from (\ref{3.66}), (\ref{3.67}) and (\ref{3.68}),
\be\label{3.69} \begin{array}{lll}
&E[\sup\limits_{s\in[t,T]}|Y_s^{0,n}-Y_s^{2,n}|^2|{\mathcal{F}}_t^{W^0}]\leq C_{\delta,\varepsilon}A_t^2\\
& +\widehat{C}_{\delta, \gamma}\varepsilon (\frac{2}{2-p})^\frac{2}{p}E[\sup\limits_{s\in[t,T]}|Y_s^{0,n}-Y_s^{2,n}|^2|{\mathcal{F}}_t^{W^0}],\ t\in[t_0,T],
\end{array}\ee
and since $\widehat{C}_{\delta, \gamma}\varepsilon (\frac{2}{2-p})^{\frac{2}{p}}<1$, we get
\be\label{3.70}E[\sup\limits_{s\in[t,T]}|Y_s^{0,n}-Y_s^{2,n}|^2|{\mathcal{F}}_t^{W^0}]\leq C_{\delta,\varepsilon}A^2_t ,\ t\in[t_0,T],\ n\geq 1,\ \mbox{P-a.s.}
\ee
Consequently, from (\ref{3.66}),
\be\label{3.71} E[\int_t^T|Z_r^{0,n}-\widehat{Z}_r^{1,n}|^2dr|\mathcal{F}_t^{W^0}]\leq C_{\delta,\varepsilon}A^2_t,\ t\in[t_0,T].\ee
\endpf
We now can prove Proposition \ref{p3.2}.

\textbf{Proof (of Proposition \ref{p3.2}).}

1) We begin with proving Assertion i). For this we note that, for all $s\in[t_0,T],\ \ n\geq 1,$
 \be\label{3.72}\begin{array}{lll}
&\widetilde{Y}_s^{1,n}-Y_s^{0,n}\leq\widehat{Y}_s^{1,n}-Y_s^{0,n} \    \ \ \  &(\mbox{Lemma}\ \ref{lemma3.4}) \\
& =A_s+(Y_s^{2,n}-Y_s^{0,n}) \ \    \ \ &(\mbox{Definition of}\ Y^{2,n})\\
&\leq A_s+|Y_s^{2,n}-Y_s^{0,n}|\\
&\leq A_s+C_\delta A_s \ \  \   \ &(\mbox{Lemma}\ \ref{Lemma3.6})\\
&\leq C_\delta(|t_0-t_1|+\sup\limits_{r\in[t_0,s]}|X_r^0-\widetilde{X}_r^1|) \ \  \hskip3cm \ \ &(\mbox{Definition of}\ A ).
\end{array}\ee
The same argument, slightly adapted, allows to show
\be\label{3.73}Y^{0,n}_s-\widetilde{Y}_s^{1,n}\leq C_\delta(|t_0-t_1|+\sup\limits_{r\in[t_0,s]}|X_r^0-\widetilde{X}_r^1|) \ , \ s\in[t_0,T],\ \ n\geq1.\ee

\noindent Thus, it only remains to prove the estimate ii) for $\widetilde{Z}^{1,n}-Z^{0,n}$, when
$$\varphi(t,x)=\varphi \in {\mathbb R},\ \ (t,x)\in[0,T]\times {\mathbb R}^d.$$
2) For this we observe that from It\^o's formula applied to $(\widetilde{Y}_t^{1,n}-Y_t^{0,n})^2$\ it follows
$$\label{3.74}\begin{array}{lll}
&  |\widetilde{Y}_s^{1,n}-Y_s^{0,n}|^2+E[\int_s^T|\widetilde{Z}_r^{1,n}-Z_r^{0,n}|^2dr|\mathcal{F}_s^{W^0}] \\
 =& E[|\Phi(\widetilde{X}_T^1)-\Phi(X_T^0)|^2|\mathcal{F}_s^{W^0}]\\
& +2E[\int_s^T(\widetilde{Y}_r^{1,n}-Y_r^{0,n})(\frac{1}{\dot{\tau}}f(\tau^{-1}(r),\widetilde{X}_r^1,\widetilde{Y}_r^{1,n},
\sqrt{\dot{\tau}}Z_r^{1,n},u_r^0)\\
&\ \mbox{}\hskip1cm-f(r,X_r^0,Y_r^{0,n},Z_r^{0,n},u_r^0))dr|\mathcal{F}_s^{W^0}]\\
&  +2\frac{n}{\dot{\tau}}E[\int_s^T(\widetilde{Y}_r^{1,n}-Y_r^{0,n})((\widetilde{Y}_r^{1,n}-\varphi)^-
-(Y_r^{0,n}-\varphi)^-)dr|\mathcal{F}_s^{W^0}]\ \mbox{} \hskip2.0cm (\leq 0)&\\
&+2n(\frac{1}{\dot{\tau}}-1)E[\int_s^T(\widetilde{Y}_r^{1,n}-Y_r^{0,n})(Y_r^{0,n}-\varphi)^-dr|\mathcal{F}_s^{W^0}]\\
\leq &  C_\delta A_s^2\ \mbox{} \hskip9.38cm     (\mbox{Lemma} \ \ref{lemma3.1})&\\
&  +C_\delta E[\int_s^T A_r(|t_0-t_1|+|\widetilde{X}_r^1-X_r^0|+|\widetilde{Y}_r^{1,n}-Y_r^{0,n}|+|\widetilde{Z}_r^{1,n}-Z_r^{0,n}|\\
&\ \mbox{}\hskip1cm+|t_0-t_1||Z_r^{0,n}|)dr
|\mathcal{F}_s^{W^0}]\ \mbox{} \hskip1.29cm(\mbox{Lemma} \ \ref{lemma3.2}\ \mbox{and i) of Proposition}\ \ref{p3.2})&\\
&  +C_\delta |t_0-t_1|(E[A_T^2|\mathcal{F}_s^{W^0}])^{\frac{1}{2}}(E[(n\int_s^T(Y_r^{0,n}-\varphi)^-dr)^2|\mathcal{F}_s^{W^0}])
^{\frac{1}{2}} \ (\mbox{Proposition} \ \ref{p3.2} \mbox{-i)}).&
\end{array}$$
Thus, again from Proposition \ref{p3.2}-i) and Lemma \ \ref{lemma3.1},
\be\label{3.74}\begin{aligned}
&  |\widetilde{Y}_s^{1,n}-Y_s^{0,n}|^2+E[\int_s^T|\widetilde{Z}_r^{1,n}-Z_r^{0,n}|^2dr|\mathcal{F}_s^{W^0}] \\
&\leq  C_\delta A_s^2+\frac{1}{2}E[\int_s^T|\widetilde{Z}_r^{1,n}-Z_r^{0,n}|^2dr|\mathcal{F}_s^{W^0}]+C_\delta|t_0-t_1|A_s
(E[\int_s^T|Z_r^{0,n}|^2dr|\mathcal{F}_s^{W^0}])^\frac{1}{2}\\
&+C_\delta|t_0-t_1|A_s(E[(n\int_s^T(Y_r^{0,n}-\varphi)^-dr)^2|\mathcal{F}_s^{W^0}])^\frac{1}{2},\ s\in [t_0, T],\ n\geq 1.
\end{aligned}\ee

\noindent Note that, due to Lemma \ref{lemma3.3}, we have $$E[\int_s^T|Z_r^{0,n}|^2dr+(n\int_s^T(Y_r^{0,n}-\varphi)^-dr)^2|\mathcal{F}_s^{W^0}]\leq C,\ s\in [t_0, T].$$
Consequently, P-a.s., for all $n\geq 1,\ s\in [t_0, T]$,
\be\label{3.75}E[\int_s^T|\widetilde{Z}_r^{1,n}-Z_r^{0,n}|^2dr|\mathcal{F}_s^{W^0}]\leq C_{\delta} A_s^2\leq C_{\delta}(|t_0-t_1|^2+\sup\limits_{r\in[t_0,s]}|\widetilde{X}_r^1-X_r^0|^2).\ee

\noindent On the other hand, recalling that $\varphi$\ is constant, from It\^o's formula, Lemma \ \ref{lemma3.1} and Proposition \ref{p3.2}-i)\ we deduce
\be\label{3.76} \begin{aligned}
& |\widetilde{Y}_s^{1,n}-Y_s^{0,n}|^2+\int_s^T|\widetilde{Z}_r^{1,n}-Z_r^{0,n}|^2 dr\\
=& (\Phi(\widetilde{X}_T^1)-\Phi(X_T^0))^2\\
&+2\int_s^T(\widetilde{Y}_r^{1,n}-Y_r^{0,n})
(\frac{1}{\dot{\tau}}f(\tau^{-1}(r),\widetilde{X}_r^1,\widetilde{Y}_r^{1,n},\sqrt{\dot{\tau}}\widetilde{Z}_r^{1,n},u_r^0)-f(r,X_r^0,Y_r^{0,n},Z_r^{0,n},u_r^0))dr\\
&+2\frac{n}{\dot{\tau}}\int_s^T(\widetilde{Y}_r^{1,n}-{Y}_r^{0,n})((\widetilde{Y}_r^{1,n}-\varphi)^--(Y_r^{0,n}-\varphi)^-)dr  \ \ (\leq 0) \\
&+2n(\frac{1}{\dot{\tau}}-1)\int_s^T(\widetilde{Y}_r^{1,n}-{Y}_r^{0,n})(Y_r^{0,n}-\varphi)^-dr-
2\int_s^T(\widetilde{Y}_r^{1.n}-Y_r^{0,n})(\widetilde{Z}_r^{1.n}-Z_r^{0,n})dW_r^0\\
&\leq CA_T^2+C_\delta\int_s^T A_r(A_r+|t_0-t_1||Z_r^{0,n}|+|Z_r^{0,n}-\widetilde{Z}_r^{1,n}|)dr\\
&+C_\delta|t_0-t_1|A_T(n\int_s^T(Y_r^{0,n}-\varphi)^-dr)-2\int_s^T(\widetilde{Y}_r^{1,n}-Y_r^{0,n})(\widetilde{Z}_r^{1,n}-Z_r^{0,n})dW_r^0\\
&\leq C_\delta A_T^2+\frac{1}{2}\int_s^T|Z_r^{0,n}-\widetilde{Z}_r^{1,n}|^2dr+|t_0-t_1|^2\int_s^T|Z_r^{0,n}|^2dr \\
&+C_\delta|t_0-t_1|A_T(n\int_s^T(Y_r^{0,n}-\varphi)^-dr)-
2\int_s^T(\widetilde{Y}_r^{1,n}-Y_r^{0,n})(\widetilde{Z}_r^{1,n}-Z_r^{0,n})dW_r^0.\end{aligned}\ee
Therefore, we have\be\label{3.77} \begin{aligned}\int_s^T|\widetilde{Z}_r^{1,n}-Z_r^{0,n}|^2dr\leq &C_\delta A_T^2+2|t_0-t_1|^2(\int_s^T|Z_r^{0,n}|^2dr+ (n\int_s^T(Y_r^{0,n}-\varphi)^-dr )^2)\\
&-4\int_s^T(\widetilde{Y}_r^{1,n}-Y_r^{0,n})(\widetilde{Z}_r^{1,n}-Z_r^{0,n})dW_r^0,\ s\in [t_0, T],\ n\geq 1,\end{aligned}\ee
and, consequently, for $p\geq 1$,
\be\label{3.78} \begin{aligned}&E[(\int_s^T|\widetilde{Z}_r^{1,n}-Z_r^{0,n}|^2 dr)^{2p}|\mathcal{F}_s^{W^0}]\leq C_pC_\delta A_s^{4p}\\
&+C_p|t_0-t_1|^{4p}\big(E[(\int_s^T|Z_r^{0,n}|^2dr)^{2p}|\mathcal{F}_s^{W^0}]
+E[(n\int_s^T(Y_r^{0,n}-\varphi)^-dr)^{4p}|\mathcal{F}_s^{W^0}]\big)\\
&+C_p E[(\int_s^T|\widetilde{Y}_r^{1,n}-Y_r^{0,n}|^2|\widetilde{Z}_r^{1,n}-Z_r^{0,n}|^2dr)^p|\mathcal{F}_s^{W^0}].
\end{aligned}\ee
We recall that, due to Lemma\ \ref{lemma3.3},
$$E[(\int_s^T|Z_r^{0,n}|^2dr)^{2p}|\mathcal{F}_s^{W^0}]\leq C_p;$$
$$E[(n\int_s^T(Y_r^{0,n}-\varphi)^-dr)^{4p}|\mathcal{F}_s^{W^0}]\leq C_p,\ s\in [t_0, T],\ n\geq 1,\ \mbox{P-a.s.}$$
Thus, due to Proposition\ \ref{p3.2}-i), $|\widetilde{Y}_r^{1,n}-Y_r^{0,n}|^2\leq C_\delta A_r^2\leq C_\delta A_T^2,\ \mbox{P-a.s.}$\ Hence, we get \be\label{3.79-1}\begin{aligned}
&E[(\int_s^T|\widetilde{Y}_r^{1,n}-Y_r^{0,n}|^2|\widetilde{Z}_r^{1,n}-Z_r^{0,n}|^2dr)^p|\mathcal{F}_s^{W^0}]\\
&\leq C_\delta E[A_T^{2p}(\int_s^T|\widetilde{Z}_r^{1,n}-Z_r^{0,n}|^2dr)^p|\mathcal{F}_s^{W^0}]\\
&\leq C_\delta(E[A_T^{6p}|\mathcal{F}_s^{W^0}])^\frac{1}{3}(E[(\int_s^T|\widetilde{Z}_r^{1,n}-Z_r^{0,n}|^2dr)^{\frac{3}{2}p}
|\mathcal{F}_s^{W^0}])^\frac{2}{3}\\
&\leq C_\delta A_s^{2p}(E[(\int_s^T|\widetilde{Z}_r^{1,n}-Z_r^{0,n}|^2dr)^{\frac{3}{2}p}
|\mathcal{F}_s^{W^0}])^\frac{2}{3}.
\end{aligned}\ee
Hence, from (\ref{3.78}) and (\ref{3.79-1}) it follows that
\be\label{3.79}E[(\int_s^T|\widetilde{Z}_r^{1,n}-Z_r^{0,n}|^2 dr)^{4p}|\mathcal{F}_s^{W^0}]\leq C_{\delta,p} A_s^{8p}+C_{\delta,p} A_s^{4p}(E[(\int_s^T|\widetilde{Z}_r^{1,n}-Z_r^{0,n}|^2dr)^{3p}|\mathcal{F}_s^{W^0}])^{\frac{2}{3}},\ee
from where we get the announced result for $p\geq\frac{1}{4}$, which means
\be\label{3.80}E[(\int_s^T|\widetilde{Z}_r^{1,n}-Z_r^{0,n}|^2dr)^p|\mathcal{F}_s^{W^0}]\leq C_{\delta,p}A_s^{2p},\ \mbox{P-a.s.},\ s\in[t_0,T],\ n\geq 1,\ p\geq1. \ee
\endpf

\section{The semiconcavity of the value function}
In this section we consider $V$ as value function of a stochastic control problem which cost functional is defined by a BSDE reflected at an upper barrier. Indeed, if it is reflected at a lower barrier, $V$ can, in general, not be semiconcave, let us illustrate this by an easy example.
\begin{example}
We consider the controlled system (\ref{FSDE}) endowed with RBSDE (\ref{RBSDE}) reflected at a lower barrier $\varphi$. $T> 1$. We let the dimension $m=d=1$\ and consider the case of no control ($U$ is a singleton) and with the coefficients  $b\equiv0,\ \sigma\equiv0,$\ $ f\equiv-1,\ \varphi\equiv0,$\ and $\Phi\equiv1.$\\
Then, obviously, $X_s^{t,x}=x,\ s\in[t,T],$\ and the solution of RBSDE (\ref{RBSDE}) is given by:
$$Y_s^{t,x}=(1-(T-s))^+,\ Z_s^{t,x}=0,\ K_s^{t,x}=(1-(T-t))^--(1-(T-s))^-,\ s\in[t,T].$$
Consequently,
$$V(t,x)=Y_t^{t,x}=(1-(T-t))^+,\ (t,x)\in[0,T]\times {\mathbb R}.$$
However, although the coefficients satisfy our assumptions, it can be easily seen that the function V is not semiconcave on $[0,T-\delta]\times {\mathbb R},$\ for all $0<\delta<1<T$.
\end{example}
For this reason, for $(t,x)\in[0,T]\times {\mathbb R}^d,\ W=(W_s)_{s\in[t,T]}$\ $m$-dimensional Brownian motion with $W_t=0,$\ and $u\in\mathcal{U}_{t,T}^W$,  we associate SDE (\ref{FSDE}) with the RBSDE reflected at an upper barrier $\varphi$:
\begin{equation}\label{Eq4.1}
\left\{
\begin{aligned}
dY_s^{t,x,u}=&\ -f(s,X_s^{t,x,u},Y_s^{t,x,u},Z_s^{t,x,u},u_s)ds+Z_s^{t,x,u}dW_s+dK_s^{t,x,u},\\
Y_T^{t,x,u}=&\Phi(X_T^{t,x,u}),\\
Y_s^{t,x,u}\leq  &\varphi(s,X_s^{t,x,u}),\ (Y_s^{t,x,u}-\varphi(s,X_s^{t,x,u}))dK_s^{t,x,u}=0,\ s\in[t,T].
\end{aligned}
\right.
\end{equation}
Under the assumptions H1) and H3') it has a unique solution
$$(Y^{t,x,u},Z^{t,x,u},K^{t,x,u}) \in  {\cal{S}}_{{\mathbb{F}}^W}^2(t, T)\times L_{{\mathbb{F}}^W}^{2}(t,T;{\mathbb{R}}^{m})\times A_{{\mathbb{F}}^W}^2(t, T) .$$
In order to emphasize the dependence on $W$, we also write
$$(Y^{t,x,u}(W), Z^{t,x,u}(W), K^{t,x,u}(W))=(Y^{t,x,u}, Z^{t,x,u}, K^{t,x,u}).$$

\noindent We define \be\label{4.2}V(t,x)=\inf_{u\in\mathcal{U}_{t,T}^W}Y_t^{t,x,u},\ (t,x)\in[0,T]\times {\mathbb R}^d,\ee
and we recall that $V\in C_b([0, T]\times {\mathbb R}^d)$\ is the unique (uniqueness in $C_p([0, T]\times {\mathbb R}^d)$) viscosity solution of the HJB equation with an upper obstacle
\begin{equation}\label{HJB-3}\left\{
\begin{aligned}
&\max\left\{V(t,x)-\varphi(t,x),-\frac{\partial}{\partial t}V(t,x)-\inf_{u\in U}H(t,x,V(t,x),\nabla V(t,x),D^2V(t,x),u)\right\}=0, \\
&V(T,x)=\Phi(x), \ \ \ \ (t,x)\in[0,T)\times {\mathbb R}^d.
\end{aligned}
\right.\end{equation}
The main result of this section is the following one.
\begin{theorem}\label{Th4.1}
We assume that the conditions H1), H2), H3'), H4) and H5) are satisfied, as well as H6) or H7). Then, for all $\delta>0$, there is some $C_\delta>0$\ such that, for all $(t_0,x_0),\ (t_1,x_1)\in[0,T-\delta]\times {\mathbb R}^d$, and for all $\lambda\in[0,1]$:
\be\label{4.3}\lambda V(t_1,x_1)+(1-\lambda)V(t_0,x_0)\leq V(\lambda(t_1,x_1)+(1-\lambda)(t_0,x_0))+C_\delta \lambda(1-\lambda)(|t_0-t_1|^2+|x_0-x_1|^2).\ee
\end{theorem}
As in Section 2, the proof will be based on the approximation of the reflected BSDE (\ref{Eq4.1}) by penalized BSDEs:
\begin{equation}\label{Eq10}
\left\{
\begin{aligned}
dY_s^{t,x,u;n}=&\ -[f(s,X_s^{t,x,u},Y_s^{t,x,u;n},Z_s^{t,x,u;n},u_s)-n(Y_s^{t,x,u;n}-\varphi(s,X_s^{t,x,u}))^+]ds+Z_s^{t,x,u;n}dW_s,\\
Y_T^{t,x,u;n}=&\Phi(X_T^{t,x,u}),\ s\in[t,T],\ n\geq 1.
\end{aligned}
\right.
\end{equation}
For every $n\geq 1$, BSDE (\ref{Eq10}) admits a unique solution $(Y^{t,x,u;n}, Z^{t,x,u;n})$, and we define
\be\label{valuefunction-n-1}V_n(t,x):=\inf_{u\in{\cal U}^W_{t,T}}Y^{t,x,u;n}_t,\ \ (t,x)\in[0,T]\times {\mathbb R}^d.
\ee

\noindent In analogy to Proposition \ref{p3.1} we have
\begin{proposition}\label{p4.1} Under the assumptions H1) and H3') the following assertions hold true:

\noindent{\rm i)} $Y_s^{t,x,u;n}\downarrow Y_s^{t,x,u},\ \mbox{as}\ n\rightarrow\infty,\ \mbox{P-a.s.},\ s\in[t,T],\ u\in\mathcal{U}_{t,T}^W$;\\
{\rm ii)} $E[\sup\limits_{s\in[t,T]}|Y_s^{t,x,u;n}-Y_s^{t,x,u}|^2+\int_t^T|Z_s^{t,x,u;n}-Z_s^{t,x,u}|^2ds
\\ \mbox{ }\ \mbox{ }\ +\sup\limits_{s\in[t,T]}|K_s^{t,x,u}-n\int_s^T(Y_r^{t,x,u;n}-\varphi(r,X_r^{t,x,u}))^+dr|^2] \rightarrow 0, \ \mbox{as}\ n\rightarrow\infty,\ u\in\mathcal{U}_{t,T}^W;$\\
{\rm iii)} $V_n(t,x)\downarrow V(t,x), \ \mbox{as}\ n\rightarrow\infty,\ (t,x)\in[0,T]\times {\mathbb R}^d$.
\end{proposition}
Theorem \ref{Th4.1} is an immediate consequence of the following theorem combined with assertion iii) of Proposition \ref{p4.1}.
\begin{theorem}\label{Th4.2}
Under the assumptions of Theorem \ref{Th4.1}, for all $\delta>0$, there is some $C_\delta\in {\mathbb R}$\ such that, for all $n\geq 1$, for all $(t_0,x_0),\ (t_1,x_1)\in[0,T-\delta]\times {\mathbb R}^d$\ and for all $\lambda \in (0,1)$,
\be\label{4.5}\lambda V_n(t_1,x_1)+(1-\lambda)V_n(t_0,x_0)\leq V_n(\lambda(t_1,x_1)+(1-\lambda)(t_0,x_0))+C_\delta\lambda(1-\lambda)(|t_0-t_1|^2+|x_0-x_1|^2).\ee
\end{theorem}

As in Section 2, our proof is based on the method of time change.\\
Let $\delta>0,\ (t_i,x_i)\in[0,T-\delta]\times {\mathbb R}^d,\ i=0,\ 1,$\ and $\lambda\in (0,1),$\ and let us put $(t_\lambda, x_\lambda):=\lambda(t_1,x_1)+(1-\lambda)(t_0,x_0)$. Moreover, let $W^\lambda=(W_s^{\lambda})_{s\in[t_\lambda,T]}$\ be a $m$-dimensional Brownian motion with $W_{t_\lambda}^{\lambda}=0$, and $u^\lambda\in\mathcal{U}_{t_\lambda,T}^{W^\lambda}$\ be an admissible control associated with $W^{\lambda}$.\\
Using the notations introduced in the preceding section, we put $$X^\lambda:=X^{t_\lambda,x_\lambda,u^\lambda}(W^\lambda);$$
$$(Y^\lambda, Z^\lambda, K^\lambda):=(Y^{t_\lambda,x_\lambda,u^\lambda}(W^\lambda), Z^{t_\lambda,x_\lambda,u^\lambda}(W^\lambda), K^{t_\lambda,x_\lambda,u^\lambda}(W^\lambda));$$
$$(Y^{\lambda,n}, Z^{\lambda,n}):=(Y^{t_\lambda,x_\lambda,u^\lambda;n}(W^\lambda), Z^{t_\lambda,x_\lambda,u^\lambda;n}(W^\lambda)),\ \ n\geq 1.$$

We use the method of time change again. But since we have to compare the stochastic control system with initial data $(t_\lambda, x_\lambda)$\ with those of initial data $(t_0, x_0)$\ and $(t_1, x_1)$, we have to define a separate time change, for each $i=0,\ 1$:
\be\label{4.6}\begin{array}{lll}
&\tau_i: [t_i,T]\rightarrow[t_\lambda,T],\ \  \tau_i(s)=t_\lambda+\frac{T-t_\lambda}{T-t_i}(s-t_i),\ s\in[t_i,T].\end{array}\ee
We observe that $\dot{\tau}_i(=\frac{d}{ds}\tau_i(s))=\frac{T-t_\lambda}{T-t_i},\ i=0,\ 1,$ and so
$ W_s^i:=\frac{1}{\sqrt{\dot{\tau}_i}}W_{\tau_i(s)}^\lambda,\ s\in[t_i,T]$, is a $m$-dimensional Brownian motion with $W_{t_i}^i=0$; and $u_s^i:=u_{\tau_i(s)}^\lambda,\ s\in[t_i,T],$\ defines an admissible control belonging to $\mathcal{U}_{t_i,T}^{W^i},\ i=0,\ 1.$

For $i=0, 1,$\ we consider the solution $X^i:=X^{t_i,x_i,u^i}(W^i)$ of SDE (\ref{FSDE}) governed by the Brownian motion $W^i$\ and the control $u^i$,\ \ as well as the solution $(Y^i, Z^i, K^i)=(Y^{t_i, x_i,u^i}(W^i),$ $ Z^{t_i,x_i,u^i}(W^i), K^{t_i,x_i,u^i}(W^i))$\ of the associated reflected BSDE (\ref{Eq4.1}),\ and
the solution $(Y^{i,n}, Z^{i,n})$ $=(Y^{t_i,x_i,u^i;n}(W^i),$ $ Z^{t_i,x_i,u^i;n}(W^i))$\ of the associated penalized BSDE (\ref{Eq10}).

We have to work with the triples $(X^\lambda, Y^{\lambda,n}, Z^{\lambda,n}),\ (X^i, Y^{i,n}, Z^{i,n}),\ \ i=0,1, \ \ n\geq 1.$\ However, in order to make them comparable, we need equations driven by the same Brownian motion. For this end we consider the inverse time changes:\\
\be\label{4.7}\begin{aligned}
&\tau_i^{-1}: [t_\lambda, T]\rightarrow[t_i, T],\ \  \tau_i^{-1}(s)=t_i+\frac{T-t_i}{T-t_\lambda}(s-t_\lambda),\ \ s\in[t_\lambda,T],\ i=0,\ 1,\end{aligned}\ee
and we introduce the time changed processes
\be\label{4.8}\begin{aligned}\widetilde{X}_s^i:=X^i_{\tau_i^{-1}(s)},\ \widetilde{Y}_s^{i,n}:=Y^{i, n}_{\tau_i^{-1}(s)},\ \ \widetilde{Z}_s^{i,n}:=\frac{1}{\sqrt{\dot{\tau}_i}}Z^{i,n}_{\tau_{i}^{-1}(s)},\ \ s\in[t_\lambda,T],\ \ i=0,\ 1.\end{aligned}\ee
By observing that
  \be\label{4.9}W^i_{\tau_i^{-1}(s)}=\frac{1}{\sqrt{\dot{\tau}_i}}W_s^\lambda,\ \ \ u^i_{\tau_i^{-1}(s)}=u_s^\lambda,\ \ \ s\in[t_\lambda, T],\ \ i=0,\ 1,\ee
we see that
    \begin{equation}\label{4.10}
\left\{
\begin{aligned}
d\widetilde{X}_s^i&=\ \frac{1}{\dot{\tau}_i}b(\tau_i^{-1}(s),\widetilde{X}_s^i,u_s^\lambda)ds+ \frac{1}{\sqrt{\dot{\tau_i}}}\sigma(\tau_i^{-1}(s),\widetilde{X}_s^i,u_s^\lambda)dW_s^\lambda,\ \ s\in[t_\lambda,T];\\
\widetilde{X}^i_{t_\lambda}&=x_i,
\end{aligned}
\right.
\end{equation}
and
 \begin{equation}\label{4.11}
\left\{
\begin{aligned}
d\widetilde{Y}_s^{i,n}&=   -(\frac{1}{\dot{\tau}_i}f(\tau_i^{-1}(s),\widetilde{X}_s^i,\widetilde{Y}_s^{i,n},\sqrt{\dot{\tau}_i}\widetilde{Z}_s^{i,n},
u_s^\lambda)- \frac{n}{\dot{\tau}_i}(\widetilde{Y}_s^{i,n}-\varphi(\tau_i^{-1}(s),\widetilde{X}_s^i))^+)ds
+\widetilde{Z}_s^{i,n}dW_s^\lambda,\\
\widetilde{Y}^{i,n}_T&= \Phi(\widetilde{X}_T^i),\ \ s\in[t_\lambda, T],\ i=0,\ 1.
\end{aligned}
\right.
\end{equation}
With the same, only slightly adapted arguments as those for Lemma \ref{lemma3.1} and Proposition \ref{p3.2}, we can show the following statement.
\begin{lemma}\label{lemma4.1} Let us suppose the assumptions H1), H2) and H3'). Then,\\
{\rm i)} For all $p\geq 1$ there is some constant $C_{\delta,p}$\ such that, for all $t\in[t_\lambda, T],\ \ n\geq 1,$\ P-a.s.,\\
(1) $ E[\sup\limits_{s\in[t,T]}|\widetilde{X}_s^0-\widetilde{X}_s^1|^p|\mathcal{F}_t^{W^\lambda}]\leq C_{\delta,p}(|t_0-t_1|^p+|\widetilde{X}_t^0-\widetilde{X}_t^1|^p)$;\\
(2) $ |\widetilde{Y}_t^{0,n}-\widetilde{Y}_t^{1,n}|\leq C_{\delta,p}(|t_0-t_1|+\sup\limits_{s\in[t_\lambda, t]}|\widetilde{X}_s^0-\widetilde{X}_s^1|)$;\\
{\rm ii)} If, moreover, $\varphi(t,x)\equiv\varphi\in {\mathbb R},\ (t,x)\in[0,T]\times {\mathbb R}^d$,
then, for all $p\geq 1$\ there is some constant $C_{\delta,p}$\ such that, for all $t\in[t_\lambda, T],\ \ n\geq 1,$\ P-a.s.,\\
$$E[\left(\int_t^T|\widetilde{Z}_r^{0,n}-\widetilde{Z}_r^{1,n}|^2dr\right)^p|\mathcal{F}_t^{W^\lambda}]\leq C_{\delta,p}\left(|t_0-t_1|^2+\sup\limits_{s\in[t_\lambda,t]}|\widetilde{X}_s^0-\widetilde{X}_s^1|^2\right)^p.$$
\end{lemma}

We also shall introduce the process $\widetilde{X}_s:=
\lambda\widetilde{X}_s^1+(1-\lambda)\widetilde{X}_s^0$,\ \
$\widetilde{Y}_s^n:=\lambda\widetilde{Y}_s^{1,n}+(1-\lambda)\widetilde{Y}_s^{0,n}$,\ \ and
$\widetilde{Z}_s^n:=\lambda\widetilde{Z}_s^{1,n}+(1-\lambda)\widetilde{Z}_s^{0,n},\ s\in[t_\lambda,T]$.
Recall the definition of the processes $X^\lambda,\ (Y^{\lambda,n},\ Z^{\lambda,n})$,
and
\be\label{4.12}A_t:=\sup_{s\in[t_\lambda,t]}(|t_0-t_1|+|\widetilde{X}_s^1-\widetilde{X}_s^0|),\ \ \
t\in[t_\lambda,T],\ee
and introduce the continuous increasing process
\be\label{4.13}B_t:=\sup_{s\in[t_\lambda,t]}|\widetilde{X}_s-X_s^\lambda|,\ \ \
t\in[t_\lambda,T],\ee we have
\bp\label{p4.2} Under the assumption of Theorem \ref{Th4.1} there is some
$C_\delta\in \mathbb{R} $\ only depending on $\delta>0$ and on the
bounds and the Lipschitz constants of $\sigma,\ b,\ f,\ \Phi,\ \varphi$,\
$\nabla_{(t,x)}\sigma$\ and $\nabla_{(t,x)}b$, such that
$$\widetilde{Y}_t^n\leqslant
Y_t^{\lambda,n}+C_\delta(B_t+\lambda(1-\lambda)A_t^2),\ t\in
[t_\lambda,T],\ n\geqslant1,\ \mbox{P-a.s.}$$\ep

Before proving Proposition \ref{p4.2} let
us show that Theorem \ref{Th4.2} holds true.\\
\textbf{Proof (of Theorem \ref{Th4.2})}.\\ We recall that $\delta>0$, and $(t_0,x_0),\ (t_1,x_1)\in
[0,T-\delta]\times \mathbb{R}^d$\ are arbitrarily chosen, and
$(t_\lambda,x_\lambda)=\lambda(t_1,x_1)+(1-\lambda)(t_0,x_0)$. For
an arbitrary $\lambda\in (0,1)$, $n\geqslant 1$, we choose $\varepsilon>0$\ small enough and we let $u^\lambda \in {\cal{U}}_{t_\lambda,T}^{W^\lambda}$\ be such that\\
\be\label{4.14}V_n(t_\lambda,x_\lambda)=\inf_{u\in
{\cal{U}}_{t_\lambda,T}^{W^\lambda}}Y_{t_\lambda}^{t_\lambda,x_\lambda,u;n}\geqslant
Y_{t_\lambda}^{t_\lambda,x_\lambda,u^\lambda;n}-\varepsilon=Y_{t_\lambda}^{\lambda,n}-\varepsilon.\ee
As $V_n(t_i,x_i)\leqslant
Y_{t_i}^{i,n}=\widetilde{Y}_{t_\lambda}^{i,n},\ i=0,\ 1$, we have from
Proposition \ref{p4.2} (note that: $B_{t_\lambda}=0,$\ and $A_{t_\lambda}=|t_0-t_1|+|x_0-x_1|$),
\be\label{4.15}\begin{array}{lll}
&\lambda V_n(t_1,x_1)+(1-\lambda)V_n(t_0,x_0)\leqslant
\lambda\widetilde{Y}_{t_\lambda}^{1,n}+(1-\lambda)\widetilde{Y}_{t_\lambda}^{0,n}=\widetilde{Y}_{t_\lambda}^n\\
&\leqslant
Y_{t_\lambda}^{\lambda,n}+C_\delta\lambda(1-\lambda)(|t_0-t_1|^2+|x_0-x_1|^2)\\
&\leqslant
V_n(t_\lambda,x_\lambda)+\varepsilon+C_\delta\lambda(1-\lambda)(|t_0-t_1|^2+|x_0-x_1|^2).\end{array}\ee
Finally, from the arbitrariness of $\varepsilon>0$,
\be\label{4.16}\lambda
V_n(t_1,x_1)+(1-\lambda)V_n(t_0,x_0)-V_n(t_\lambda,x_\lambda)\leqslant C_\delta\lambda(1-\lambda)(|t_0-t_1|^2+|x_0-x_1|^2).\ee
Note that $C_\delta$ does neither depend on
$\lambda,\ (t_0,x_0)$\ and $(t_1,x_1)$, nor on $n\geqslant1$.\endpf

The proof of Proposition \ref{p4.2} is split into a sequel of lemmas. The following lemma will be crucial for our computations.
\bl\label {lemma4.2}  For all $p\geqslant1$ there is some $C_{p,\delta}\in
\mathbb{R}$ only depending on $\delta,p$ and on the bounds and the
Lipschitz constants of $\sigma$ and $b$, such that, $t\in[t_\lambda,T], \mbox{P-a.s.},$
\be\label{4.17}E[\sup_{s\in[t,T]}|\widetilde{X}_s-X_s^\lambda|^p|{\cal{F}}_t^{W^\lambda}]\leqslant C_p|\widetilde{X}_t-X_t^\lambda|^p+
C_{p,\delta}(\lambda(1-\lambda))^p(|t_0-t_1|^2+|\widetilde{X}_t^1-\widetilde{X}_t^0|^2)^p.\ \ee\el
\noindent \textbf{Proof}. For $s\in[t_\lambda, T]$, we have to estimate the equation
\be\label{4.18}\begin{aligned}
d(\widetilde{X}_s-X_s^\lambda)&=\bigg(\frac{\lambda}{\dot{\tau}_1}b\big(\tau_1^{-1}(s),\widetilde{X}_s^1,
u_s^\lambda\big)+\frac{1-\lambda}{\dot{\tau}_0}b\big(\tau_0^{-1}(s),\widetilde{X}_s^0,u_s^\lambda\big)-
b(s,X_s^\lambda,u_s^\lambda)\bigg)ds\\
& +\bigg(\frac{\lambda}{\sqrt{\dot{\tau}_1}}\sigma\big(\tau_1^{-1}(s),\widetilde{X}_s^1,u_s^\lambda\big)
+\frac{1-\lambda}{\sqrt{\dot{\tau}_0}}\sigma\big(\tau_0^{-1}(s),\widetilde{X}_s^0,u_s^\lambda\big)-
\sigma(s,X_s^\lambda,u_s^\lambda)\bigg)dW_s^\lambda,\\
\widetilde{X}_{t_\lambda}-X_{t_\lambda}^\lambda&=\lambda\widetilde{X}_{t_\lambda}^1+(1-\lambda)
\widetilde{X}_{t_\lambda}^0-X_{t_\lambda}^\lambda=0.\end{aligned}\ee
For this let us begin with\\
1) Estimating
$\big|\big(\frac{\lambda}{\sqrt{\dot{\tau}_1}}\sigma\big(\tau_1^{-1}(s),\widetilde{X}_s^1,u_s^\lambda\big)
+\frac{1-\lambda}{\sqrt{\dot{\tau}_0}}\sigma\big(\tau_0^{-1}(s),\widetilde{X}_s^0,u_s^\lambda\big)\big)
-\sigma(s,X_s^\lambda,u_s^\lambda)\big|$.\\
From a straight-forward computation we get
$$\lambda|1-\frac{1}{\sqrt{\dot{\tau}_1}}|\leqslant\frac{1}{2\delta}\lambda(1-\lambda)|t_0-t_1|; \ \ (1-\lambda)|1-\frac{1}{\sqrt{\dot{\tau}_0}}|\leqslant\frac{1}{2\delta}\lambda(1-\lambda)|t_0-t_1|;$$
and
$$|\lambda(1-\frac{1}{\sqrt{\dot{\tau}_1}})+(1-\lambda)(1-\frac{1}{\sqrt{\dot{\tau}_0}})|\leqslant
\frac{1}{\delta^2}\lambda(1-\lambda)|t_1-t_0|^2.$$
We also observe that $|\tau_1^{-1}(s)-\tau_0^{-1}(s)|\leqslant|t_1-t_0|,\ \
s\in[t_\lambda, T]$.\\
Consequently,
\be\label{4.19}\begin{array}{rcl}
& &\big|\lambda(1-\frac{1}{\sqrt{\dot{\tau}_1}})\sigma\big(\tau_1^{-1}(s),\widetilde{X}_s^1,u_s^\lambda\big)
+(1-\lambda)(1-\frac{1}{\sqrt{\dot{\tau}_0}})\sigma\big(\tau_0^{-1}(s),\widetilde{X}_s^0,u_s^\lambda\big)\big|\\
& \leqslant&\lambda\big|1-\frac{1}{\sqrt{\dot{\tau}_1}}\big|\big|\sigma\big(\tau_1^{-1}(s),\widetilde{X}_s^1,
u_s^\lambda\big)-\sigma\big(\tau_0^{-1}(s),\widetilde{X}_s^0,u_s^\lambda\big)\big|+C\big|
\lambda(1-\frac{1}{\sqrt{\dot{\tau}_1}})+(1-\lambda)(1-\frac{1}{\sqrt{\dot{\tau}_0}})\big|\\
&\leqslant&
C_\delta\lambda(1-\lambda)(|t_0-t_1|^2+|\widetilde{X}_s^0-\widetilde{X}_s^1|^2),\ \ s\in[t_\lambda,T].\end{array}\ee
Also remark that, thanks to assumption H5) the functions
$\sigma(\cdot,\cdot,u),\ (-\sigma)(\cdot,\cdot,u),\ b(\cdot,\cdot,u),$\ $(-b)(\cdot,\cdot,u)$\
are semiconcave, uniformly with respect to $u\in  {U}$. Thus, from the latter estimate
\be\label{4.20}\begin{array}{rcl}
& &\big|\frac{\lambda}{\sqrt{\dot{\tau}_1}}\sigma\big(\tau_1^{-1}(s),\widetilde{X}_s^1,u_s^\lambda\big)
+\frac{1-\lambda}{\sqrt{\dot{\tau}_0}}\sigma\big(\tau_0^{-1}(s),\widetilde{X}_s^0,u_s^\lambda\big)
-\sigma(s,X_s^\lambda,u_s^\lambda)\big|\\
&\leqslant&\big|\lambda\sigma\big(\tau_1^{-1}(s),\widetilde{X}_s^1,u_s^\lambda\big)
+(1-\lambda)\sigma\big(\tau_0^{-1}(s),\widetilde{X}_s^0,u_s^\lambda\big)
-\sigma(s,X_s^\lambda,u_s^\lambda)\big|\\
& &+C_\delta\lambda(1-\lambda)(|t_1-t_0|^2+|\widetilde{X}_s^1-\widetilde{X}_s^0|^2)\\
&\leqslant&\big|\lambda\sigma\big(\tau_1^{-1}(s),\widetilde{X}_s^1,u_s^\lambda\big)
+(1-\lambda)\sigma\big(\tau_0^{-1}(s),\widetilde{X}_s^0,u_s^\lambda\big)-
\sigma\big(\lambda\tau_1^{-1}(s)+(1-\lambda)\tau_0^{-1}(s),\widetilde{X}_s,u_s^\lambda\big)\big|\\
& &+\big|\sigma\big(\lambda\tau_1^{-1}(s)+(1-\lambda)\tau_0^{-1}(s),\widetilde{X}_s,u_s^\lambda\big)
-\sigma(s,X_s^\lambda,u_s^\lambda)\big|\\
& &+C_\delta\lambda(1-\lambda)(|t_1-t_0|^2
+|\widetilde{X}_s^1-\widetilde{X}_s^0|^2)\\
&\leqslant&
C_\delta\lambda(1-\lambda)(|t_1-t_0|^2+|\widetilde{X}_s^1-\widetilde{X}_s^0|^2)\\
& &
+C(|\lambda\tau_1^{-1}(s)+(1-\lambda)\tau_0^{-1}(s)-s|+|\widetilde{X}_s-X_s^\lambda|),\ \ s\in[t_\lambda,T].\end{array}\ee
However,
$\lambda\tau_1^{-1}(s)+(1-\lambda)\tau_0^{-1}(s)-s\equiv0,s\in[t_\lambda,T]$,
so that
\be\label{4.21}\begin{array}{rcl}
& &\big|\frac{\lambda}{\sqrt{\dot{\tau}_1}}\sigma\big(\tau_1^{-1}(s),
\widetilde{X}_s^1,u_s^\lambda\big)+\frac{1-\lambda}{\sqrt{\dot{\tau}_0}}\sigma\big(\tau_0^{-1}(s),
\widetilde{X}_s^0,u_s^\lambda\big)-\sigma(s,X_s^\lambda,u_s^\lambda)\big|\\
& &\leqslant
C_\delta\lambda(1-\lambda)(|t_1-t_0|^2+|\widetilde{X}_s^1-\widetilde{X}_s^0|^2)
+C|\widetilde{X}_s-X_s^\lambda|,\ \ s\in[t_\lambda,T].\end{array}\ee
2) By using now
$$\lambda(1-\frac{1}{\dot{\tau}_1})=\lambda(1-\lambda)\frac{t_1-t_0}{T-t_\lambda},\ \
(1-\lambda)(1-\frac{1}{\dot{\tau}_0})=\lambda(1-\lambda)\frac{t_0-t_1}{T-t_\lambda},$$
we get similarly to (\ref{4.21}):\\
\be\label{4.22}\begin{array}{rcl}
& &\big|\frac{\lambda}{\dot{\tau}_1}b\big(\tau_1^{-1}(s),\widetilde{X}_s^1,u_s^\lambda\big)
+\frac{1-\lambda}{{\dot{\tau}_0}}b\big(\tau_0^{-1}(s),\widetilde{X}_s^0,u_s^\lambda\big)
-b(s,X_s^\lambda,u_s^\lambda)\big|\\
&\leqslant & C_\delta\lambda(1-\lambda)(|t_0-t_1|^2+|\widetilde{X}_s^1-\widetilde{X}_s^0|^2)
+C|\widetilde{X}_s-X_s^\lambda|,\ \ s\in[t_\lambda,T].\end{array}\ee
From (\ref{4.21}), (\ref{4.22}), Lemma \ref{lemma4.1} and standard SDE estimates we then get the
wished result.\endpf

Now we have still to prepare the proof of Proposition \ref{p4.2}. For this we recall that\be\label{4.33}\begin{array}{rcl}
d\widetilde{Y}_s^n&=&-\Big\{\frac{\lambda}{\dot{\tau}_1}f(\tau_1^{-1}(s),
\widetilde{X}_s^1,\widetilde{Y}_s^{1,n},\sqrt{\dot{\tau}_1}\widetilde{Z}_s^{1,n},u_s^\lambda)
+\frac{1-\lambda}{\dot{\tau}_0}f(\tau_0^{-1}(s),\widetilde{X}_s^0,\widetilde{Y}_s^{0,n},
\sqrt{\dot{\tau}_0}\widetilde{Z}_s^{0,n},u_s^\lambda)\\
& &-\Big(\lambda\frac{n}{\dot{\tau}_1}\big(\widetilde{Y}_s^{1,n}-\varphi(\tau_1^{-1}(s),\widetilde{X}_s^1)\big)^+
+(1-\lambda)\frac{n}{\dot{\tau}_0}\big(\widetilde{Y}_s^{0,n}
-\varphi(\tau_0^{-1}(s),\widetilde{X}_s^0)\big)^+\Big)\Big\}ds\\
& &+\widetilde{Z}_s^ndW_s^\lambda,\\
\widetilde{Y}_T^n&=&\lambda\Phi(\widetilde{X}_T^1)+(1-\lambda)\Phi(\widetilde{X}_T^0),\end{array}\ee
and we compare this equation with the BSDE
\be\label{4.34}\begin{array}{rcl}
d\widehat{Y}_s^n&=&-\Big(f(s,X_s^\lambda,\widehat{Y}_s^n-CB_s-C_\delta\lambda(1-\lambda)A_s^2,
\widehat{Z}_s^n,u_s^\lambda)+C(CB_s+C_\delta\lambda(1-\lambda)A_s^2)\\
& &+C_\delta^0\lambda(1-\lambda)\big(|t_0-t_1|^2(1+|\widetilde{Z}_s^{0,n}|^2)
+|\widetilde{Z}_s^{1,n}-\widetilde{Z}_s^{0,n}|^2\big)\\
& &-n\big(\widehat{Y}_s^n-\varphi(s,X_s^\lambda)-CB_s-C_\delta\lambda(1-\lambda)A_s^2\big)^+\Big)ds
+\widehat{Z}_s^ndW_s^\lambda,\\
\widehat{Y}_T^n&=&\Phi(X_T^\lambda)+CB_T+C_\delta\lambda(1-\lambda)A_T^2,\end{array}\ee
where $C_\delta^0=0$, if $f$ is independent of $z$.

\br We point out that, due to Lemma \ref{lemma3.3},
$$E[(\int_s^T|\widetilde{Z}_r^{i,n}|^2dr)^p|{\cal{F}}_s^{W^\lambda}]\leqslant
C_{\delta, p},\ s\in [t_\lambda, T],\ \ p\geqslant1.$$ This shows that above BSDE (\ref{4.34}) is well-posed.\er

\bl\label{lemma100} Under the assumptions for Theorem \ref{Th4.1} we have
$$\widetilde{Y}_t^n\leqslant\widehat{Y}_t^n,\ \ t\in[t_\lambda, T],\ n\geqslant1,\ \mbox{P-a.s.}$$
\el
\noindent \textbf{Proof}. The proof is based on the comparison theorem (Lemma \ref{comparisonbsde} in Section 4). We prepare for the application of this comparison theorem by the following three steps.\\
\textbf{Step 1}.\\
Using that $-a^+-b^+\leqslant -(a+b)^+,\  a,\ b\in\mathbb{R}$, we have
\be\label{4.23}\begin{array}{rcl}
& &-\lambda\frac{n}{\dot{\tau}_1}\big(\widetilde{Y}_t^{1,n}-\varphi(\tau_1^{-1}(t),\widetilde{X}_t^1)\big)^+
-(1-\lambda)\frac{n}{\dot{\tau}_0}\big(\widetilde{Y}_t^{0,n}-\varphi(\tau_0^{-1}(t),\widetilde{X}_t^0)\big)^+\\
&\leqslant&-n\Big(\frac{\lambda}{\dot{\tau}_1}\widetilde{Y}_t^{1,n}
+\frac{1-\lambda}{\dot{\tau}_0}\widetilde{Y}_t^{0,n}-
\big(\frac{\lambda}{\dot{\tau}_1}\varphi(\tau_1^{-1}(t),\widetilde{X}_t^1)
+\frac{1-\lambda}{\dot{\tau}_0}\varphi(\tau_0^{-1}(t),\widetilde{X}_t^0)\big)\Big)^+\\
&=&-n\big\{\widetilde{Y}_t^n-\big(\lambda(1-\frac{1}{\dot{\tau}_1})\widetilde{Y}_t^{1,n}
+(1-\lambda)(1-\frac{1}{\dot{\tau}_0})\widetilde{Y}_t^{0, n}\big)
-\big(\lambda\varphi(\tau_1^{-1}(t),\widetilde{X}_t^1)+(1-\lambda)\varphi(\tau_0^{-1}(t),\widetilde{X}_t^0)\big)\\
& &+\lambda(1-\frac{1}{\dot{\tau}_1})\varphi(\tau_1^{-1}(t),\widetilde{X}_t^1)+(1-\lambda)
(1-\frac{1}{\dot{\tau}_0})\varphi(\tau_0^{-1}(t),\widetilde{X}_t^0)\big\}^+\\
&=&-n\big\{\widetilde{Y}_t^n-\frac{\lambda(1-\lambda)}{T-t_\lambda}(t_1-t_0)
(\widetilde{Y}_t^{1,n}-\widetilde{Y}_t^{0,n})-\big(\lambda\varphi(\tau_1^{-1}(t),\widetilde{X}_t^1)
+(1-\lambda)\varphi(\tau_0^{-1}(t),\widetilde{X}_t^0)\big)\\
& &+\frac{\lambda(1-\lambda)}{T-t_\lambda}(t_1-t_0)\big(\varphi(\tau_1^{-1}(t),\widetilde{X}_t^1)-
\varphi(\tau_0^{-1}(t),\widetilde{X}_t^0)\big)\big\}^+\\
&\leqslant&-n\big\{\widetilde{Y}_t^n-\big(\lambda\varphi(\tau_1^{-1}(t),\widetilde{X}_t^1)
+(1-\lambda)\varphi(\tau_0^{-1}(t),\widetilde{X}_t^0)\big)-C_\delta\lambda(1-\lambda)|t_1-t_0|A_t\big\}^+, \end{array}\ee
where Lemma \ref{lemma4.1}\ was applied for the latter inequality.\\
Hence, from the semiconcavity of $\varphi$, and since $\lambda\tau_1^{-1}(t)+(1-\lambda)\tau_0^{-1}(t)=t,$
\be\label{4.24}\begin{array}{rcl}
& &-\lambda\frac{n}{\dot{\tau}_1}\big(\widetilde{Y}_t^{1,n}-\varphi(\tau_1^{-1}(t),\widetilde{X}_t^1)\big)^+
-(1-\lambda)\frac{n}{\dot{\tau}_0}\big(\widetilde{Y}_t^{0,n}-\varphi(\tau_0^{-1}(t),\widetilde{X}_t^0)\big)^+\\
&\leqslant &-n\big(\widetilde{Y}_t^n-\varphi(t,\widetilde{X}_t)-C_\delta\lambda(1-\lambda)A_t^2\big)^+\\
&\leqslant &-n\big(\widetilde{Y}_t^n-\varphi(t,X_t^\lambda)-CB_t-C_\delta\lambda(1-\lambda)A_t^2\big)^+,\ \ t\in[t_\lambda,T],\ \ n\geq 1,\end{array}\ee
where
$B_t:=\sup_{s\in[t_\lambda,t]}|\widetilde{X}_s-X_s^\lambda|$.\\
We recall that, from Lemma \ref{lemma4.2}
\be\label{4.25}E[B_T^p|{\cal{F}}_t^{W^\lambda}]\leqslant
C_{p}B_t^p+C_{p,\delta}(\lambda(1-\lambda))^p(|t_1-t_0|^2+|\widetilde{X}_t^1-\widetilde{X}_t^0|^2)^p,\
t\in[t_\lambda,T],\ p\geqslant1,\ \mbox{P-a.s.} \ee
Hence,
\be\label{4.26}\begin{array}{rcl}
& &-\lambda\frac{n}{\dot{\tau}_1}\big(\widetilde{Y}_t^{1,n}-\varphi(\tau_1^{-1}(t),\widetilde{X}_t^1)\big)^+
-(1-\lambda)\frac{n}{\dot{\tau}_0}\big(\widetilde{Y}_t^{0,n}-\varphi(\tau_0^{-1}(t),\widetilde{X}_t^0)\big)^+\\
&\leqslant &-n\big(\widetilde{Y}_t^n-\varphi(t,X_t^\lambda)-CB_t-C_\delta\lambda(1-\lambda)A_t^2\big)^+,\
t\in[t_\lambda,T],\ n\geqslant1.\end{array}\ee
Note that, if $\varphi$ is a constant independent of $(t,x)$, then
\be\label{4.27}\begin{array}{rcl}
& &-\lambda\frac{n}{\dot{\tau}_1}(\widetilde{Y}_t^{1,n}-\varphi)^+-(1-\lambda)\frac{n}{\dot{\tau}_0}(\widetilde{Y}_t^{0,n}-\varphi)^+\\
&\leqslant &-n\big(\widetilde{Y}_t^n-\varphi-C_\delta\lambda(1-\lambda)A_t^2\big)^+,\ t\in[t_\lambda,T],\ n\geqslant1.\end{array}\ee

\noindent\textbf{Step 2}.\\ From the semiconcavity of $f$ and standard arguments similar to those used in Step 1 we obtain,
\be\label{4.28}\begin{aligned}
&  \frac{\lambda}{\dot{\tau}_1}f(\tau_1^{-1}(s),\widetilde{X}_s^1,\widetilde{Y}_s^{1,n},
\sqrt{\dot{\tau}_1}\widetilde{Z}_s^{1,n},u_s^\lambda)+\frac{1-\lambda}{\dot{\tau}_0}f(\tau_0^{-1}(s),
\widetilde{X}_s^0,\widetilde{Y}_s^{0,n},\sqrt{\dot{\tau}_0}\widetilde{Z}_s^{0,n},u_s^\lambda)\\
= &\lambda
f(\tau_1^{-1}(s),\widetilde{X}_s^1,\widetilde{Y}_s^{1,n},\sqrt{\dot{\tau}_1}\widetilde{Z}_s^{1,n},u_s^\lambda)
+(1-\lambda)f(\tau_0^{-1}(s),\widetilde{X}_s^0,\widetilde{Y}_s^{0,n},\sqrt{\dot{\tau}_0}\widetilde{Z}_s^{0,n},
u_s^\lambda)\\
&  -\lambda(1-\lambda)\frac{t_1-t_0}{T-t_\lambda}\{f(\tau_1^{-1}(s),\widetilde{X}_s^1,\widetilde{Y}_s^{1,n},
\sqrt{\dot{\tau}_1}\widetilde{Z}_s^{1,n},u_s^\lambda)-f(\tau_0^{-1}(s),\widetilde{X}_s^0,\widetilde{Y}_s^{0,n},
\sqrt{\dot{\tau}_0}\widetilde{Z}_s^{0,n},u_s^\lambda)\}\\
\leqslant &
f(s,\widetilde{X}_s,\widetilde{Y}_s^n,\lambda\sqrt{\dot{\tau}_1}\widetilde{Z}_s^{1,n}
+(1-\lambda)\sqrt{\dot{\tau}_0}\widetilde{Z}_s^{0,n},u_s^\lambda) +C_\delta\lambda(1-\lambda)(|t_1-t_0|^2+|\widetilde{X}_s^1-\widetilde{X}_s^0|^2\\
&  +|\widetilde{Y}_s^{1,n}-\widetilde{Y}_s^{0,n}|^2+|\widetilde{Z}_s^{1,n}-\widetilde{Z}_s^{0,n}|^2
+|t_1-t_0|^2|\widetilde{Z}_s^{0,n}|^2),\ \ s\in[t_\lambda,T].\end{aligned}\ee
Since, on the other hand,
\be\label{4.29}\begin{array}{rcl}
& &|\lambda\sqrt{\dot{\tau}_1}\widetilde{Z}_s^{1,n}+(1-\lambda)\sqrt{\dot{\tau}_0}\widetilde{Z}_s^{0,n}-\widetilde{Z}_s^n|\\
&=&|\lambda(1-\sqrt{\dot{\tau}_1})\widetilde{Z}_s^{1,n}+(1-\lambda)(1-\sqrt{\dot{\tau}_0})\widetilde{Z}_s^{0,n}|\\
&\leqslant&\lambda|1-\sqrt{\dot{\tau}_1}||\widetilde{Z}_s^{1,n}-\widetilde{Z}_s^{0,n}|
+|\lambda(1-\sqrt{\dot{\tau}_1})+(1-\lambda)(1-\sqrt{\dot{\tau}_0})||\widetilde{Z}_s^{0,n}|\\
&\leqslant&
C_\delta\lambda(1-\lambda)(|t_1-t_0||\widetilde{Z}_s^{1,n}-\widetilde{Z}_s^{0,n}|+|t_1-t_0|^2|\widetilde{Z}_s^{0,n}|)
\end{array}\ee
(see the proof of Lemma \ref{lemma4.2}), we have
\be\label{4.30}\begin{array}{rcl}
& &\frac{\lambda}{\dot{\tau}_1}f(\tau_1^{-1}(s),\widetilde{X}_s^1,\widetilde{Y}_s^{1,n},
\sqrt{\dot{\tau}_1}\widetilde{Z}_s^{1,n},u_s^\lambda)+\frac{1-\lambda}{\dot{\tau}_0}f(\tau_0^{-1}(s),
\widetilde{X}_s^0,\widetilde{Y}_s^{0,n},\sqrt{\dot{\tau}_0}\widetilde{Z}_s^{0,n},u_s^\lambda)\\
&\leqslant &
f(s,\widetilde{X}_s,\widetilde{Y}_s^n,\widetilde{Z}_s^n,u_s^\lambda)
+C_\delta\lambda(1-\lambda)\big(|t_0-t_1|^2(1+|\widetilde{Z}_s^{0,n}|^2)+|\widetilde{X}_s^1-\widetilde{X}_s^0|^2\\
& &+|\widetilde{Y}_s^{1,n}-\widetilde{Y}_s^{0,n}|^2+|\widetilde{Z}_s^{1,n}-\widetilde{Z}_s^{0,n}|^2\big),\ \ s\in[t_\lambda,T],\ \ n\geqslant1.\end{array}\ee
Thus, from Lemmas \ref{lemma4.1} and \ref{lemma4.2}\\
\be\label{4.31}\begin{array}{rcl}
& &\frac{\lambda}{\dot{\tau}_1}f(\tau_1^{-1}(s),\widetilde{X}_s^1,\widetilde{Y}_s^{1,n},
\sqrt{\dot{\tau}_1}\widetilde{Z}_s^{1,n},u_s^\lambda)+\frac{1-\lambda}{\dot{\tau}_0}
f(\tau_0^{-1}(s),\widetilde{X}_s^0,\widetilde{Y}_s^{0,n},\sqrt{\dot{\tau}_0}\widetilde{Z}_s^{0,n},u_s^\lambda)\\
&\leqslant&
f(s,X_s^\lambda,\widetilde{Y}_s^n-CB_s-C_\delta\lambda(1-\lambda)A_s^2,\widetilde{Z}_s^n,u_s^\lambda)+C'B_s+C'_\delta
\lambda(1-\lambda)A_s^2\\
& &+C_\delta^0\lambda(1-\lambda)\big(|t_0-t_1|^2(1+|\widetilde{Z}_s^{0,n}|^2)
+|\widetilde{Z}_s^{1,n}-\widetilde{Z}_s^{0,n}|^2\big),\ \ s\in[t_\lambda,T],\ \ n\geqslant1.\end{array}\ee
Remark also that if $f$\ does not depend on $z$, the constant $C_\delta^0$\ in (\ref{4.31}) can be chosen to be zero.

\noindent\textbf{Step 3}. \\
We also note that, thanks to the semiconcavity and the Lipschitz condition on $\Phi$,\  \be\label{100}\lambda\Phi(\widetilde{X}_T^1)+(1-\lambda)\Phi(\widetilde{X}_T^0)\leqslant\Phi(X_T^\lambda)
+C_\delta\lambda(1-\lambda)A_T^2+CB_T.\ee
The above three steps allow to conclude. Indeed, taking into account (\ref{4.26}), (\ref{4.31}), and (\ref{100}), it follows from the comparison
theorem-Lemma \ref{comparisonbsde} in Section 4 that:\be\label{4.35}\widetilde{Y}_t^n\leqslant\widehat{Y}_t^n,\ t\in[t_\lambda,T],\ n\geqslant1.\ee
The proof is complete.\endpf.

Let us now introduce the process
$$\overline{Y}_t^n=\widehat{Y}_t^n-CB_t-C_\delta\lambda(1-\lambda)A_t^2,\ \ t\in[t_\lambda,T].$$
Then, for
$D_t:=CB_t+C_\delta\lambda(1-\lambda)A_t^2,\ \ t\in[t_\lambda,T]$, we have
\be\left\{\label{4.36}\begin{array}{rcl}
d\overline{Y}_s^n&=&-\big\{f(s,X_s^\lambda,\overline{Y}_s^n,\widehat{Z}_s^n,u_s^\lambda)+CD_s
+C_\delta^0\lambda(1-\lambda)\big(|t_0-t_1|^2(1+|\widetilde{Z}_s^{0,n}|^2)
+|\widetilde{Z}_s^{1,n}-\widetilde{Z}_s^{0,n}|^2\big)\\
& &-n\big(\overline{Y}_s^n-\varphi(s,X_s^\lambda)\big)^+\big\}ds-dD_s+\widehat{Z}_s^ndW_s^\lambda,\ \ s\in[t_\lambda,T],\\
\overline{Y}_T^n&=&\Phi(X_T^\lambda).\end{array}\right.\ee
Recalling that\be\left\{\label{4.37}\begin{array}{rcl}dY_s^{\lambda,n}&=&-\big\{f(s,X_s^{\lambda},
Y_s^{\lambda,n},Z_s^{\lambda,n},u_s^\lambda)-n\big(Y_s^{\lambda,n}-\varphi(s,X_s^\lambda)\big)^+\big\}ds\\
& &+Z_s^{\lambda,n}dW_s^\lambda,\ \ s\in[t_\lambda,T],\\
Y_T^{\lambda,n}&=&\Phi(X_T^\lambda).\end{array}\right.\ee
We can establish the following statement.
\bl\label{lemma4.3} Under the assumptions of Theorem \ref{Th4.1}
$$E[\sup_{s\in[t,T]}|\overline{Y}_s^n-Y_s^{\lambda,n}|^2+\int_t^T|\widehat{Z}_s^n-Z_s^{\lambda,n}|^2ds|
{\cal{F}}_t^{W^\lambda}]\leqslant
C_\delta D_t^2,\ \ t\in[t_\lambda,T].$$\el

\noindent \textbf{Proof}. Taking into account that
$$-(\overline{Y}_s^n-Y_s^{\lambda,n})\big((\overline{Y}_s^n-\varphi(s,X_s^\lambda))^+
-(Y_s^{\lambda,n}-\varphi(s,X_s^\lambda))^+\big)\leqslant0,\ s\in [t_\lambda,T],$$
we see that, for $\gamma>0$,
\be\label{4.38}\begin{aligned}
& e^{\gamma
t}(\overline{Y}_t^n-Y_t^{\lambda,n})^2+E[\int_t^Te^{\gamma
s}\big(\gamma|\overline{Y}_s^n-Y_s^{\lambda,n}|^2+|\widehat{Z}_s^n-Z_s^{\lambda,n}|^2\big)ds|{\cal{F}}_t^{W^\lambda}]\\
 \leqslant& 2E[\int_t^Te^{\gamma
s}(\overline{Y}_s^n-Y_s^{\lambda,n})\big(f(s,X_s^\lambda,\overline{Y}_s^n,\widehat{Z}_s^n,u_s^\lambda)
-f(s,X_s^{\lambda},Y_s^{\lambda,n},Z_s^{\lambda,n},u_s^\lambda)\big)ds|{\cal{F}}_t^{W^\lambda}]\\
  &+2E[\int_t^Te^{\gamma
s}(\overline{Y}_s^n-Y_s^{\lambda,n})dD_s|{\cal{F}}_t^{W^\lambda}]\\
  &+2E[\int_t^Te^{\gamma
s}(\overline{Y}_s^n-Y_s^{\lambda,n})\Big(CD_s+C_\delta^0\lambda(1-\lambda)
\big(|t_0-t_1|^2(1+|\widetilde{Z}_s^{0,n}|^2)+|\widetilde{Z}_s^{1,n}-\widetilde{Z}_s^{0,n}|^2\big)\Big)ds
|{\cal{F}}_t^{W^\lambda}]\\
 \leqslant& CE[\int_t^Te^{\gamma
s}|\overline{Y}_s^n-Y_s^{\lambda,n}|^2ds|{\cal{F}}_t^{W^\lambda}]+\frac{1}{2}E[\int_t^Te^{\gamma
s}|\widehat{Z}_s^n-Z_s^{\lambda,n}|^2ds|{\cal{F}}_t^{W^\lambda}]+C_\gamma
E[D_T^2|{\cal{F}}_t^{W^\lambda}]\\
  &+C_{\gamma,\delta}
E[\sup_{s\in[t,T]}|\overline{Y}_s^n-{Y}_s^{\lambda,n}|
\Big(\lambda(1-\lambda)\big(|t_0-t_1|^2(1+\int_t^T|\widetilde{Z}_s^{0,n}|^2ds)
+\int_t^T|\widetilde{Z}_s^{1,n}-\widetilde{Z}_s^{0,n}|^2ds)\\
  &+D_T\Big)|{\cal{F}}_t^{W^\lambda}].\end{aligned}\ee
Let
$$\widetilde{D}_{t,T}:=\lambda(1-\lambda)\big(|t_0-t_1|^2(1+\int_t^T|\widetilde{Z}_s^{0,n}|^2ds)
+\int_t^T|\widetilde{Z}_s^{1,n}-\widetilde{Z}_s^{0,n}|^2ds)+D_T.$$
Recall that $C_\delta^0=0$, if $f$ does not depend on $z$. If $f$
depends on $z$, we have thanks to assumption H7) that $\varphi$
is constant (see Lemma \ref{lemma4.1}-ii):
$$E[(\int_t^T|\widetilde{Z}_s^{0,n}-\widetilde{Z}_s^{1,n}|^2ds)^p|{\cal{F}}_t^{W^\lambda}]\leqslant
C_{\delta,p}A_t^{2p},\ p\geqslant1,\ \mbox{P-a.s.}$$
On the other hand, from Lemma \ref{lemma3.3}-ii) we know that
$$E[(\int_t^T|{Z}_s^{0,n}|^2ds)^p|{\cal{F}}_t^{W^0}]\leqslant
C_p,\ t\in[t_0,T],\ n\geqslant1,\ p\geqslant1,\ \mbox{P-a.s.},$$
and, hence,
$$E[(\int_t^T|\widetilde{Z}_s^{0,n}|^2ds)^p|{\cal{F}}_t^{W^\lambda}]\leqslant
C_{\delta,p},\ t\in[t_\lambda,T],\ n\geqslant1,\ p\geqslant1,\ \mbox{P-a.s.}$$
Consequently, considering that from the Lemmas \ref{lemma4.1}-i) and \ref{lemma4.2} it follows that  $E[D_T^2|{\cal{F}}_t^{W^\lambda}]\leqslant C_\delta D_t^2$, we get, for $\gamma\geqslant C+1$,
\be\label{4.39}\begin{array}{rcl}
& &(\overline{Y}_t^n-Y_t^{\lambda,n})^2+E[\int_t^T|\widehat{Z}_s^n-Z_s^{\lambda,n}|^2ds|{\cal{F}}_t^{W^\lambda}]\\
&\leqslant & C_\delta
D_t^2+C_\gamma E[\sup_{s\in[t,T]}|\overline{Y}_s^n-Y_s^{\lambda,n}|\widetilde{D}_{t,T}|{\cal{F}}_t^{W^\lambda}],\ \ t\in[t_\lambda,T].\end{array}\ee
Finally, applying the argument used for (\ref{3.70-1}) in the proof of Lemma \ref{lemma3.5} (or, (\ref{3.66}) in the proof of Lemma \ref{Lemma3.6}) it
follows that
$$E[\sup_{s\in[t,T]}|\overline{Y}_s^n-Y_s^{\lambda,n}|^2+\int_t^T|\widehat{Z}_s^n-Z_s^{\lambda,n}|^2ds|{\cal{F}}_t^{W^\lambda}]\leqslant
C_\delta D_t^2,\ \ t\in[t_\lambda,T],\ n\geq 1.$$
The proof is complete.\endpf
Lemmas \ref{lemma100} and \ref{lemma4.3} allow to give the proof of Proposition \ref{p4.2}.\\
\noindent\textbf{Proof (of Proposition \ref{p4.2})}.\\
From Lemmas \ref{lemma100} and \ref{lemma4.3} we can conclude that,\be\label{4.40}\begin{array}{rcl}
\widetilde{Y}_t^n& (=&\lambda\widetilde{Y}_t^{1,n}+(1-\lambda)\widetilde{Y}_t^{0,n})\leqslant
\widehat{Y}_t^n\\
& &=\overline{Y}^n_t+D_t
 \leqslant Y_t^{\lambda,n}+C_\delta
D_t\\
& &=Y_t^{\lambda,n}+C_\delta(CB_t+C_\delta\lambda(1-\lambda)A_t^2),\ t\in[t_\lambda,T],\ n\geqslant1,\ \mbox{P-a.s.}\end{array}\ee
Thus, the proof is complete now. \endpf.
\section{ {\large Appendix}}
\subsection{{\large BSDEs}}

  The objective of this section is to recall some basic results concerning backward and reflected backward SDEs, which are frequently used in our paper. Let $(\Omega, {\cal F}, P)$\ be a compact probability space endowed with a
$d$-dimensional Brownian motion and let $T>0$\ be a finite time horizon. By
${\mathbb{F}}^{W}=\{{\mathcal{F}^{W}}_s,\ 0\leq s \leq T\}$\ we denote the
natural filtration generated by the Brownian motion $W$\ and
augmented by all P-null sets, i.e.,
$${\mathcal{F}}_s=\sigma\{W_r, r\leq s\}\vee {\mathcal{N}}_P,\ \  s\in [0, T]. $$
Here $ {\cal{N}}_P$ is the set of all P-null sets.

 A measurable function $g:
\Omega\times[0,T]\times {\mathbb{R}} \times {\mathbb{R}}^{d}
\rightarrow {\mathbb{R}} $\ satisfies that $(g(t, y,
z))_{t\in [0, T]}$ is ${\mathbb{F}}$-progressively measurable for
all $(y,z)$ in ${\mathbb{R}} \times {\mathbb{R}}^{d}$, and also the following standard assumptions:
 \vskip0.2cm

(A1) There is some real $C\ge 0$  such that, P-a.s., for all $t\in
[0, T],\
y_{1}, y_{2}\in {\mathbb{R}},\ z_{1}, z_{2}\in {\mathbb{R}}^d,\\
\mbox{ }\hskip4cm   |g(t, y_{1}, z_{1}) - g(t, y_{2}, z_{2})|\leq
C(|y_{1}-y_{2}| + |z_{1}-z_{2}|).$
 \vskip0.2cm

(A2) $g(\cdot,0,0)\in L^{2}_{\mathbb{F}}(0,T;{\mathbb{R}})$. \vskip0.2cm

 The following result on BSDEs is well-known now, for its proof the reader is referred
 to the pioneering paper by
 Pardoux and Peng~\cite{pp}.
 \bl Let the function $g$ satisfy the assumptions (A1) and (A2). Then, for any random variable $\xi\in L^2(\O, {\cal{F}}_T,$ $P),$ the
BSDE
 \be Y_t = \xi + \int_t^Tg(s,Y_s,Z_s)ds - \int^T_tZ_s\,
dW_s,\q 0\le t\le T, \label{BSDE2.1} \ee
 has a unique adapted solution
$$(Y_t, Z_t)_{t\in [0, T]}\in {\cal{S}}_{\mathbb{F}^{W}}^2(0, T)\times
L_{\mathbb{F}^{W}}^{2}(0,T;{\mathbb{R}}^{d}). $$\el

\noindent Besides the above existence and uniqueness result we also recall the important comparison
theorem for BSDEs (see, e.g., Theorem 2.2 in El Karoui, Peng,
Quenez~\cite{epq} or Proposition 2.4 in Peng~\cite{Pe2}).

\bl (Comparison Theorem)\label{comparisonbsde} Given two coefficients $g_1$ and $g_2$
satisfying (A1) and (A2) and two terminal values $ \xi_1,\ \xi_2 \in
L^{2}(\Omega, {\cal{F}}_{T}, P)$, we denote by $(Y^1,Z^1)$\ and
$(Y^2,Z^2)$\ the solution of the BSDE with the data $(\xi_1,g_1 )$\
and $(\xi_2,g_2 )$, respectively. Then we have:

{\rm (i) }(Monotonicity) If  $ \xi_1 \geq \xi_2$  and $ g_1 \geq
g_2, \ a.s.$, then $Y^1_t\geq Y^2_t$, for all $t\in [0, T]$, a.s.

{\rm (ii)}(Strict Monotonicity) If, in addition to {\rm (i)}, we
also assume that $P\{\xi_1 > \xi_2\}> 0$, then $P\{Y^1_t>
Y^2_t\}>0,$ for all $\ 0 \leq t \leq T,$\ and in particular, $
Y^1_0> Y^2_0.$ \el

\subsection{{\large Reflected BSDEs}}
After the above very short recall on BSDEs let us come now to reflected BSDEs (RBSDES). Here we only introduce RBSDEs with lower barriers; the results on RBSDEs with upper barriers are symmetric. An RBSDE is connected with a terminal value $\xi \in
L^{2}(\Omega,{\cal{F}}_{T}, P)$, a generator $g$\ and a ``barrier"
process $\{S_t\}_{0\leq t \leq
   T}$. We shall make the following condition on the barrier process:

\bs
 (A3)\ $\{S_t\}_{0\leq t \leq T}$ is a continuous process such that $\{S_t\}_{0\leq t \leq T}\in {\cal{S}}_{\mathbb{F}^{W}}^2(0, T)$.
\bs

A solution of an RBSDE is a triple $(Y, Z, K)$\ which is
${\mathbb{F}}$-progressively measurable processes, take its values
in $\mathbb{R}\times\mathbb{R}^d\times\mathbb{R}_+$, and
   satisfy the following conditions

\medskip

 {\rm (i)} $Y \in {\cal{S}}^2(0, T; {\mathbb{R}}), \, Z \in
 {\cal{H}}^{2}(0,T;{\mathbb{R}}^{d})$\ and $K_{T} \in L^{2}(\Omega,{\cal{F}}_{T}, P)$;
\be \mbox{\rm (ii)}  \ Y_t = \xi + \int_t^Tg(s,Y_s,Z_s)ds + K_{T} -
K_{t} - \int^T_tZ_sdW_s,\quad t\in [0,T];\qquad\qquad\qquad\
\label{RBSDE2.2}\ee

{\rm (iii)} $Y_t \geq S_t$,\ a.s., for any $ t\in [0,T];$
\vskip0.2cm

{\rm (iv)} $\{K_{t}\}$ is continuous and increasing, $K_{0}=0$ and $\displaystyle
\int_0^T(Y_t - S_t)dK_{t}=0.$
\medskip

The following two lemmas can be referred to Theorem 5.2 and Theorem 4.1 in El Karoui, Kapoudjian, Pardoux, Peng and
Quenez \cite{ekppq}, respectively.

 \bl Assume that $g$ satisfies (A1) and (A2), $ \xi$\ belongs to $L^{2}(\Omega, {\cal{F}}_{T}, P)$, $\{S_t\}_{0\leq t \leq T}$
satisfies (A3), and $S_T \leq \xi\ \ a.s.$ Then RBSDE~(\ref{RBSDE2.2})
has a unique solution $(Y, Z, K)\in {\cal{S}}_{{\mathbb{F}}^W}^2(0, T)\times L_{{\mathbb{F}}^W}^{2}(0,T;{\mathbb{R}}^{d})\times A_{{\mathbb{F}}^W}^2(0, T).$\el

\br For simplicity, a given triple $(\xi, g, S)$ is said to satisfy
the Standard Assumptions if the coefficient $g$ satisfies (A1) and
(A2), the terminal condition $ \xi$\ belongs to $ L^{2}(\Omega,
{\cal{F}}_{T}, P)$, the barrier process $S$\ satisfies (A3) and
$S_T \leq \xi, \ \mbox{a.s.}$ \er

\bl (Comparison Theorem)\label{comparisonrbsde} Assume that two triples $(\xi_1, g_1,
S^1)$ and $(\xi_2, g_2, S^2)$\ satisfy the Standard Assumptions, and one of the both generators $g_1$\ and $ g_2$\ to
be Lipschitz. Furthermore, we make the following assumptions:
$$
  \begin{array}{ll}
{\rm(i)}&\xi_1 \leq \xi_2,\ \ a.s.;\\
{\rm(ii)}&g_1(t,y,z) \leq g_2(t,y,z),\ a.s., \hbox{ \it for } (t,y,z)\in [0,T]\times {\mathbb{R}}\times {\mathbb{R}}^{d};\\
{\rm(iii)}& S_t^1 \leq S^2_t,\ \ a.s., \hbox{ \it for } t\in [0,T]. \\
 \end{array}
  $$
  Let $(Y^1,Z^1, K^1)$ and $(Y^2, Z^2, K^2)$ be solutions of RBSDEs~(\ref{RBSDE2.2}) with data $(\xi_1, g_1,
  S^1)$ and $(\xi_2, g_2, S^2),$ respectively.  Then $Y^1_{t} \leq Y^2_{t},\ a.s., $ for $t\in [0,T].$\el

\bl  Let $(Y,Z,K)$ be the solution of the above RBSDE~(\ref{RBSDE2.2})
with data $(\xi, g, S)$\ satisfying the above Standard Assumptions. Then
there exists a constant $C$\ such that
$$
 E [\sup_{t\leq s\leq T }
|Y_s|^2+\int_t^T|Z_s|^2ds+|K_T-K_t|^2 |{{\cal {F}}_t} ]\leq CE
[\xi^2+\left(\int_t^T g(s,0,0)ds\right)^2+\sup_{t\leq s\leq T} S_s^2
|{{\cal {F}}_t} ]. $$ The constant $C$\ depends only on the
Lipschitz constant of $g$.\el

\bl Let $(\xi,g,S)$ and $(\xi^{\prime},g^{\prime},S^{\prime})$ be
two triples satisfying the above Standard Assumptions. $(Y,Z,K)$\ and $(Y^{\prime},Z^{\prime},K^{\prime})$\ are the
solutions of RBSDE~(\ref{RBSDE2.2}) with the data $(\xi,g,S)$\ and
$(\xi^{\prime},g^{\prime},S^{\prime})$, respectively. We define
$$\Delta\xi=\xi-\xi^{\prime},\qquad \Delta g=g-g^{\prime},\qquad \Delta S=S-S^{\prime};$$
$$\Delta Y=Y-Y^{\prime},\qquad \Delta Z=Z-Z^{\prime},\qquad \Delta K=K-K^{\prime}.$$
Then there
exists a constant $C$ such that, $$\aligned & E[ \sup_{t\leq s\leq T}|\Delta Y_s|^2
+\int_t^T|\Delta Z_s|^2ds +|\Delta K_T-\Delta K_t|^2|{{\cal {F}}_t} ]\\
&\leq C E[|\Delta\xi|^2+\left(\int_t^T |\Delta
g(s,Y_s,Z_s)|ds\right)^2|{{\cal {F}}_t} ]+C\left( E[\sup_{t\leq
s\leq T}|\Delta S_s|^2|{{\cal {F}}_t}
]\right)^{1/2}\Psi_{t,T}^{1/2},
\endaligned$$
{\it where} $$\aligned \Psi_{t,T} &=
E[|\xi|^2+\left(\int_t^T|g(s,0,0)|ds
\right)^2 +\sup_{t\leq s\leq T} |S_s|^2\\
&\quad +|\xi^{\prime}|^2+\left(\int_t^T|g^{\prime}(s,0,0)|ds
\right)^2 +\sup_{t\leq s\leq T} |S^{\prime}_s|^2|{{\cal {F}}_t} ].
\endaligned$$
The constant $C$\ depends only on the Lipschitz constant of $g'$.
\el

The Lemmas 4.5 and 4.6 refer to the Propositions 3.5 and 3.6 in~\cite{ekppq}, and their generalizations can be consulted in~\cite{WY} (the Propositions 2.1 and 2.2 therein), respectively.

\br For the Markovian case where the barrier process is a deterministic function of the solution of the associated forward equation, Lemma 4.6 has been considerably improved. Indeed, Proposition 6.1 in \cite{bl1} shows that $Y$ is Lipschitz with respect to the possibly random initial condition of the driving forward SDE which solution governs the RBSDE as well as its barrier.\er

\section*{Acknowledgments} Rainer Buckdahn and Juan Li thank the Department of Applied
Mathematics of The Hong Kong Polytechnic University, P. R. China, for its hospitality.

\end{document}